\documentclass{article}

\usepackage{amssymb}

\usepackage{graphicx}

\newtheorem{thm}{Theorem}[section]
\newtheorem{prop}[thm]{Proposition}
\newtheorem{lem}[thm]{Lemma}
\newtheorem{cor}[thm]{Corollary}

\newcommand{\N}{\mathbb{N}}  % The natural numbers.
\newcommand{\R}{\mathbb{R}}  % The real numbers.
\begin{document}
	
\title{A finite atlas for solution manifolds of differential systems with discrete state-dependent delays} 
	
\author{Hans-Otto Walther}
	
%\address{Mathematisches Institut, Universit\"{a}t %Gie{\ss}en,
%Arndtstr. 2, D 35392 Gie{\ss}en, Germany. E-mail {\tt
%Hans-Otto.Walther@math.uni-giessen.de}}
	
%\thanks{Mathematisches Institut, Universit\"{a}t Gie{\ss}en,
%Arndtstr. 2, D 35392 Gie{\ss}en, Germany, E-mail: {\tt
%Hans-Otto.Walther@math.uni-giessen.de}}
	
%\setcounter{equation}{-1}	

%\begin{abstract}

%\end{abstract}
	
%\begin{keyword}
%Delay differential equation \sep state-dependent delay %\sep solution manifold
%\MSC[]{34K05, 34K43}
%\end{keyword}
	
\maketitle
	
%% Should the address above be replaced by the 
%% home address ?

\begin{abstract}
	
Let $r>0, n\in\N, {\bf k}\in\N$.
Consider the delay differential equation
$$
x'(t)=g(x(t-d_1(Lx_t)),\ldots,x(t-d_{{\bf k}}(Lx_t)))
$$
for $g:(\R^n)^{{\bf k}}\supset V\to\R^n$ continuously differentiable, $L$ a continuous linear map from $C([-r,0],\R^n)$
into a finite-dimensional vectorspace $F$, each  $d_k:F\supset W\to[0,r]$, $k=1,\ldots,{\bf k}$, continuously differentiable, and $x_t(s)=x(t+s)$. The solutions define a semiflow of continuously differentiable solution operators on the submanifold
$X_f\subset C^1([-r,0],\R^n)$ which is given by the compatibility condition $\phi'(0)=f(\phi)$ with 
$$f
(\phi)=g(\phi(-d_1(L\phi)),\ldots,\phi(-d_{{\bf k}}(L\phi))).
$$ 
We prove that $X_f$ has a finite atlas of at most $2^{{\bf k}}$ manifold charts, whose domains are 
{\it almost graphs over} $X_0$. The size of the atlas depends solely on the zerosets of the delay functions $d_k$.
\end{abstract}

\bigskip

\noindent
Key words: Delay differential equation,  state-dependent delay, solution manifold

\medskip
\noindent
MSC: 34K05, 34K43

\section{Introduction}

For given  
$r>0$ and $n\in\N$ let $C_n=C([-r,0],\R^n)$ and $C^1_n=C^1([-r,0],\R^n)$ denote the Banach spaces of continuous and continuously differentiable maps $\phi:[-r,0]\to\R^n$, respectively, with the norms given by 
$$
|\phi|=\max_{-r\le t\le0}|\phi(t)|\quad\mbox{and}\quad|\phi|_1=|\phi|+|\phi'|.
$$
By a solution manifold we mean a set of the form
$$
X_f=\{\phi\in U:\phi'(0)=f(\phi)\}
$$
with a map $f:C^1_n\supset U\to\R^n$. Suppose $U\subset C^1_n$ is open, $f$ is continuously differentiable, $X_f\neq\emptyset$,

\medskip

(e) {\it each derivative $Df(\phi):C^1_n\to\R^n$, $\phi\in U$, has a linear extension $D_ef(\phi):C_n\to\R^n$, and the map
$$
U\times C_n\ni(\phi,\chi)\mapsto D_ef(\phi)\chi\in\R^n
$$
is continuous.}

\medskip

Then $X_f$ is a continuously submanifold of codimension $n$ in $C^1_n$ \cite{W1,HKWW}. It is on such manifolds that differential equations with state-dependent delay, like for example
\begin{equation}
x'(t)=g(x(t),x(t-d)),\quad d=\rho(x(t)),
\end{equation}
with continuously differentiable functions $g:\R^2\supset V\to\R$ and $\rho:\R\to[0,r]$, define a semiflow of continuously differentiable solution operators. Equations of the form (1) have been studied in numerous papers, see 
\cite{N,KS,M-PN1,M-PN2,M-PN5,M-PN6,KA,Ke1,St,KW}.
 
\medskip

The extension property (e) is a relative of the notion of being {\it almost Fr\'echet differentiable} which was introduced by Mallet-Paret, Nussbaum, and Paraskevopoulos in \cite{M-PNP}. 

\medskip

The present paper is concerned with the nature of solution manifolds which are associated with systems of differential equations with discrete state-dependent delays. For these solution manifolds we construct a finite atlas of manifold charts, whose size does not depend on the number of equations in the system but on the the zerosets of the delays. 

\medskip

It is convenient to introduce the following terminology, for a continuously differentiable submanifold $X$ of a Banach space $E$, an open subset ${\mathcal O}\subset E$ with $X\cap{\mathcal O}\neq\emptyset$, and a closed subspace $H$ with a closed complementary space. We say $X\cap{\mathcal O}$ is a graph (over $H$) if there are a closed complementary space $Q$ for $H$ and a continuously differentiable map $\gamma:H\supset dom\to Q$ with
$$
X\cap{\mathcal O}=\{\zeta+\gamma(\zeta)\in E:\zeta\in dom\}.
$$
$X\cap{\mathcal O}$ is called an {\it almost graph} (over $H$)
if there is a continuously differentiable map $\alpha:H\supset dom\to E$ with
\begin{eqnarray*}
\alpha(\zeta) & = & 0\quad\mbox{on}\quad dom\,\cap(X\cap{\mathcal O}),\\
\alpha(\zeta) & \in & E\setminus H\quad\mbox{on}\quad dom\setminus (X\cap{\mathcal O}),\\
\end{eqnarray*}
so that the map 
$$
dom\ni\zeta\mapsto\zeta+\alpha(\zeta)\in E
$$ 
defines a diffeomorphism onto $X\cap{\mathcal O}$.

A diffeomorphism $A:{\mathcal O}\to E$ onto an open subset 
of $E$ is called an {\it almost graph diffeomorphism} (associated with ${\mathcal O}$, $X$, and $H$) if
$$
A(X\cap{\mathcal O})=H\cap A({\mathcal O})
$$
and
$$
(ag)\qquad A(\zeta)=\zeta\quad\mbox{on}\quad H\cap(X\cap{\mathcal O}).
$$

\medskip

The first one of the previous equations corresponds to what is required locally for $X$ to be a continuously differentiable submanifold of $E$. Property (ag) yields that $X\cap{\mathcal O}$ is an almost graph over $H$. In order to see this, set $dom=A(X\cap{\mathcal O})=H\cap A({\mathcal O})$ and define $\alpha:H\supset dom\to E$ by $\alpha(\zeta)=A^{-1}(\zeta)-\zeta$. 
Because of
$$
\zeta+\alpha(\zeta)=A^{-1}(\zeta)\quad\mbox{for}\quad\zeta\in dom=A(X\cap{\mathcal O})
$$
the map $dom\ni\zeta\mapsto\zeta+\alpha(\zeta)\in E$ defines a diffeomorphism onto $X\cap{\mathcal O}$.
For $\zeta\in dom\cap(X\cap{\mathcal O})\subset H\cap(X\cap{\mathcal O})$ property (ag) yields $A(\zeta)=\zeta$, hence $A^{-1}(\zeta)=\zeta$, and thereby, $\alpha(\zeta)=0$.
For $\zeta\in dom\setminus(X\cap{\mathcal O})$ we have $A^{-1}(\zeta)\in X\cap{\mathcal O}$ and $\zeta\in H$. Suppose
$\alpha(\zeta)\in H$. Then $A^{-1}(\zeta)=\alpha(\zeta)+\zeta\in H$. Hence $A^{-1}(\zeta)\in H\cap(X\cap{\mathcal O})$. Using property (ag) we infer
$\zeta=A(A^{-1}(\zeta))=A^{-1}(\zeta)\in X\cap{\mathcal O}$ which contradicts the choice of $\zeta$. It follows that $\alpha(\zeta)\in E\setminus H$ on $dom\setminus(X\cap{\mathcal O})$.

\medskip

As an example for an almost graph diffeomorphism in the plane $E=\R\cdot e_1\oplus\R\cdot e_2$ let us mention the map $A:{\mathcal O}\to\R\cdot e_1+(-1,1)\cdot e_2$ given by   
\begin{eqnarray*}
{\mathcal O} & = & (0,2)\cdot S^1\setminus[0,\infty)\cdot e_2,\\
A(x) & = & (s(y),r-1)),\quad x=ry,\quad y\in S^1\setminus\{e_2\},\quad0<r<2, 
\end{eqnarray*}
with the stereographic projection
$$
s:S^1\setminus\{e_2\}\ni y\mapsto\frac{y_1}{1-y_2}\in\R.
$$
The map $A$ is an almost graph diffeomorphism associated with its domain ${\mathcal O}$, $X=S^1\setminus\{e_2\}$, and $H=\R\cdot e_1$.

\medskip

The finite atlas obtained in the present paper consists of manifold charts whose domains are almost graphs over the closed subspace
$$
X_0=\{\phi\in C^1_n:\phi'(0)=0\}
$$
of codimension $n$ (a complementary subspace is spanned by the maps $\phi_{\nu}\in C^1_n$, $\nu=1,\ldots,n$,  with components given by $\phi_{\nu\nu}(t)=t$ and $\phi_{\nu\mu}(t)=0$ for $\nu\neq\mu$). Obviously, $X_0=X_f$ for the zero map $f:C^1_n\ni\phi\mapsto 0\in\R^n$.  

\medskip

The solution manifolds considered in the sequel are given by maps $f:U\to\R^n$ which represent the right hand side of systems of the form
\begin{equation}
x'(t)=g(x(t-d_1(Lx_t)),\ldots,x(t-d_{{\bf k}}(Lx_t))),
\end{equation}
with the solution segment $x_t:[-r,0]\to\R^n$, $t\in\R$, defined by $x_t(s)=x(t+s)$ provided $[t-r,t]$ belongs to the domain of the map $x$. In Eq. (2) $L$ is a continuous linear map from $C_n$  onto a finite-dimensional vectorspace $F$ (with the canonical topology), the delay functionals $d_k$, $k\in\{1,\ldots,{\bf k}\}=K$, are continuously differentiable maps from an open set $W\subset F$ into the interval $[0,r]$, and $g$ is a continuously differentiable map from an open set $V\subset(\R^n)^{{\bf k}}$ into $\R^n$. In order that Eq. (2) makes sense we assume that 

\medskip

(V) {\it there exist $\phi\in C^1_n$ with $L\phi\in W$ and
$(\phi(-d_1(L\phi)),\ldots,\phi(-d_{{\bf k}}(L\phi)))\in V$.}

\medskip

With the abbreviation 
$$
\widehat{\phi}=(\phi(-d_1(L\phi)),\ldots,\phi(-d_{{\bf k}}(L\phi)))
$$
we get that
\begin{equation}
U=\{\phi\in C^1_n:L\phi\in W\quad\mbox{and}\quad\widehat{\phi}\in V\}
\end{equation}
is non-empty, and Eq. (2) can be written as
\begin{equation}
x'(t)=f(x_t)
\end{equation}
with $f:U\to\R^n$ given by
\begin{equation}
f(\phi)=g(\widehat{\phi}).
\end{equation}
Proposition 2.1 below guarantees that $U$ is an open subset of $C^1_n$ and that $f$ is continuously differentiable with property (e). In Proposition 2.3 we obtain among others $X_f\neq\emptyset$, which in combination with the properties of $f$  implies that $X_f$  is a continuously differentiable submanifold of codimension $n$ in $C^1_n$. The proof of Proposition 2.3
relies on Lemma 2.2, which is instrumental also for the main results  
in Sections 3 and 4 below. Before describing their content we introduce some notation related to the zerosets of the delay functions $d_k$. For a subset $J\subset K$ let
$$
W_J=\{w\in W:d_j(w)=0\quad\mbox{for}\quad j\in J\quad\mbox{and}\quad d_k(w)>0\quad\mbox{for}\quad k\in K\setminus J\}.
$$
The sets $W_J$, $J\subset K$, are mutually disjoint,
and
$$
W=\cup_{J\subset K}W_J.
$$
Notice that $W_{\emptyset}$ is open in $F$ and that $W_K$ is closed in $W$.

\begin{figure}[h]
%\centering
\scalebox{.5}{\includegraphics{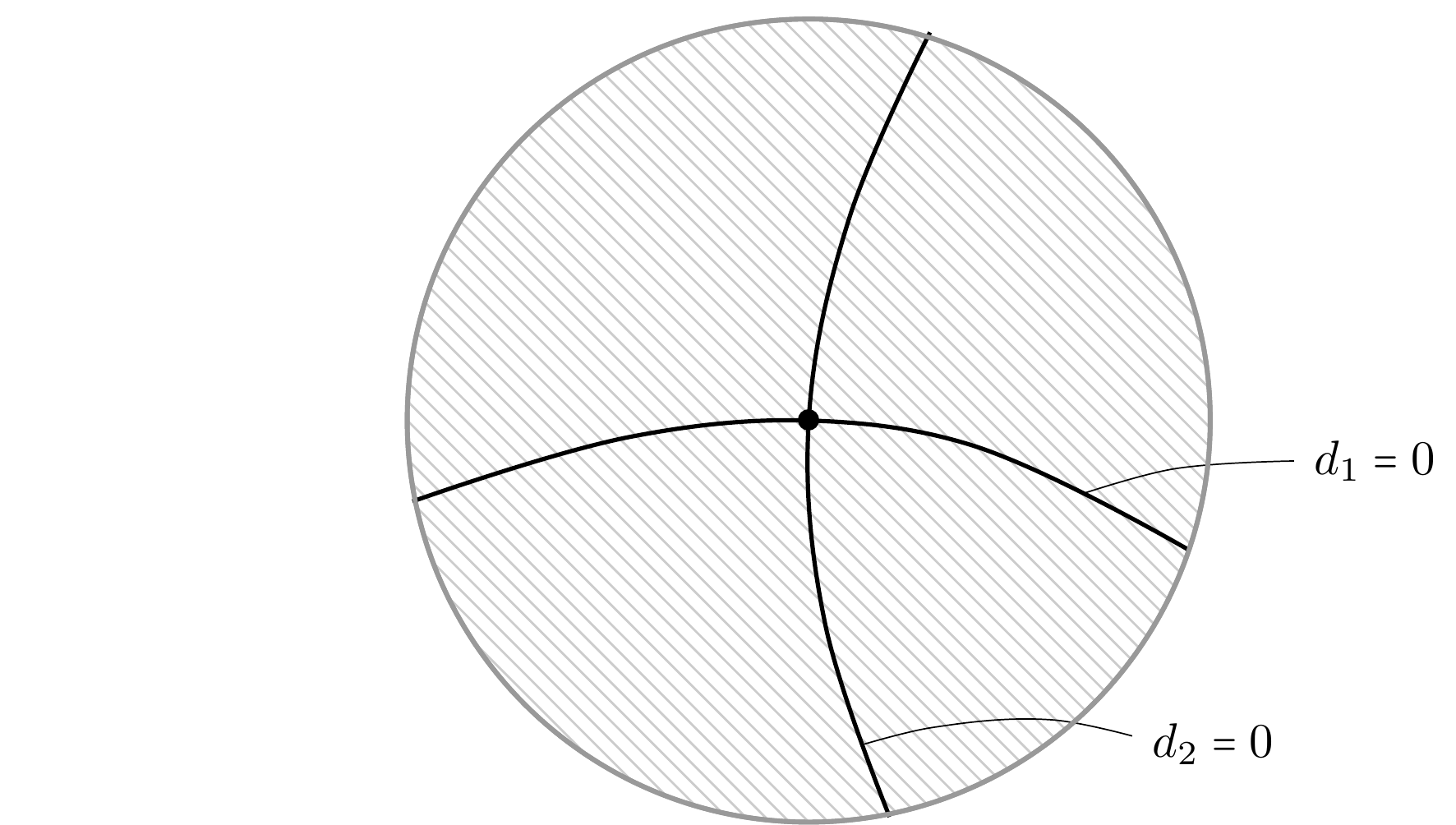}}
\caption{Decomposition of $W$ into the sets $W_{\emptyset}$ (shaded), $W_{\{1\}}$ and $W_{\{2\}}$ (black lines), and $W_{\{1,2\}}$ (the dot) in case $K=\{1,2\}$.}
%	\label{Figure 9.1}
\end{figure}

For $J\subset K$ the set 
$$
W^J=\{w\in W:d_k(w)>0\quad\mbox{for}\quad k\in K\setminus J\}
$$
is open in $F$ and contains $W_J$. Let
$$
U_J=U\cap L^{-1}(W_J),\quad X_{fJ}=X_f\cap U_J,\quad\mbox{and}\quad U^J=U\cap L^{-1}(W^J).
$$
According to the hypothesis (V) at least one of the sets  $U_J$ is non-empty. Proposition 2.3 below asserts that $X_{fJ}\neq\emptyset$ for $J\subset K$ with $U_J\neq\emptyset$. Henceforth
$$
X_f=\cup_{J\subset K:U_J\neq\emptyset}X_{fJ}
$$
is a decomposition into mutually disjoint non-empty subsets.

\begin{figure}[h]
%\centering
\scalebox{.5}{\includegraphics{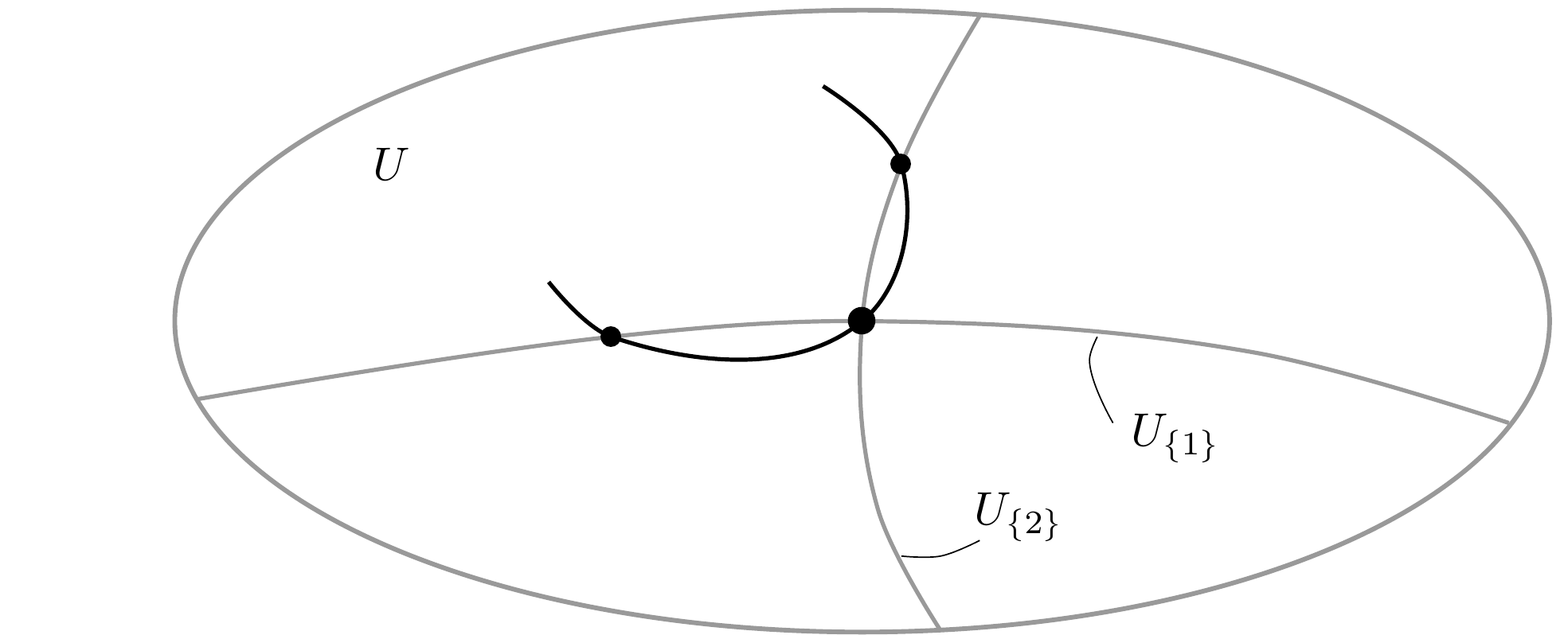}}
\caption{Decomposition of the solution manifold $X_f\subset U$ in case $K=\{1,2\}$, into $X_{f\emptyset}$ (black lines), $X_{f\{1\}}$ and $X_{f\{2\}}$ (small dots), and $X_{f\{1,2\}}$ (fat dot).}
\end{figure}

 Section 3 deals with the set $X_{fK}$ in case $U_K\neq\emptyset$. With the aid of Lemma 2.2 we find a complementary space $Q_K$ for $X_0$ so that the associated projection $C^1_n\to C^1_n$ onto $X_0$  is injective on $X_{fK}$. In Proposition 3.4 injectivity is extended to an open neighbourhood of $X_{fK}$ in $X_f$, under an {\it additional hypothesis}. Upon that we obtain in Theorem 3.5 a diffeomorphism 
 $A_K:{\mathcal O}_K\to C^1_n$ with
 $A_K(X_f\cap{\mathcal O}_K)=X_0\cap A_K({\mathcal O}_K)$ which on  $X_f\cap{\mathcal O}_K$ coincides with the said projection. The map $A_K$ is a simple almost graph diffeomorphism.
 
 \medskip
 
 Section 4 is about $X_{fJ}$ with $U_J\neq\emptyset$ for a subset $J\neq K$ of $K$ and requires more effort.
 Using Lemma 2.2 we find complementary spaces $Q_{J\chi}$ for $X_0$, for each $\chi\in X_0$ with $L\chi\in W^J$. The spaces $Q_{J\chi}$ yield a map $R^J$ from $U^J$ into $X_0$ which satisfies $R^J(\chi)=\chi$ on $X_0\cap U^J$ (like a projection) and is injective on $X_{fJ}$. The same {\it additional hypothesis} as in Section 3 allows us to extend injectivity to a neighbourhood of $X_{fJ}$ in $X_f$. The map $R^J$ defines a manifold chart on a smaller neighbourhood $N_J\subset X_f$, which is an almost graph over $X_0$ (Theorem 4.7). In Theorem 4.8 we obtain an almost graph diffeomorphism which is
 associated with an open neighbourhood ${\mathcal O}_J$ of $X_{fJ}$ in $C^1_n$, and with $X_f$ and $X_0$, and which coincides with $R^J$ on $X_f\cap{\mathcal O}_J\subset N_J$.

\medskip 

An example from \cite{W6} serves to show that unlike $X_f\cap {\mathcal O}_K$ from Section 3 the almost graph $X_f\cap{\mathcal O}_J$ from Theorem 4.8, with $J\neq K$, can in general not be written
as a graph with respect to any direct sum decomposition of $C^1_n$ into closed subspaces.

\medskip

The almost graphs $X_f\cap{\mathcal O}_J\supset X_{fJ}$, $J\subset K$ and $U_J\neq\emptyset$, are domains of manifold charts and cover $X_f$, so we obtain a finite atlas of $X_f$ with size
$$
\#\,\{J\subset K:U_J\neq\emptyset\}\le 2^{{\mathbf k}}.
$$

\medskip

As for the {\it additional hypothesis} used in Sections 3 and 4 each one of the following three properties is sufficient:

\medskip

(b) {\it $g$ is bounded},

\medskip

(d1b) {\it $d_1(w)=0$ for every $w\in W$, and for every bounded set $B\subset\R^n$ the set
$$	
\{g(v)\in\R^n:v=(v_1,\ldots,v_{\mathbf k})\in  V\subset(\R^n)^{{\mathbf k}}\,\,\mbox{and}\,\,v_1\in B\}
$$
is bounded},

\medskip

(J) $U_J=U$ for some $J\subset K$. 

\medskip

Incidentally, notice that in case of property (J) for some subset $J\subset K$  the finite atlas consists of a single chart, and the solution manifold is an almost graph over $X_0$. If for example no delay vanishes anywhere,
$$
d_k(w)>0\quad\mbox{for all}\quad w\in W\quad\mbox{and}\quad k\in K,
$$
then we have $U_{\emptyset}=U$.

\medskip

Theorems 3.5, 4.7 and 4.8 apply to Eq. (1), which has the form (2) with $n=1,r>0,K=\{1,2\},F=\R$, and
$$
L\phi=\phi(0),\quad d_1(w)=0,\quad d_2(w)=\rho(w).
$$ 
In case $\rho(\xi)>0$ everywhere the hypothesis $U_{\{1\}}=U$ is satisfied. See work of Stumpf \cite{St} for a particular example.  In case
the function $\rho$ has zeros and $g$ is bounded condition  (b) is obvious. If $\rho$ has zeros and Eq. (1) has the form
$$
x'(t)=g(x(t),x(t-d))=a\,x(t)+h(x(t-d)),\quad d=\rho(x(t)),
$$ 
with a real parameter $a$ and a bounded continuously differentiable function $h$ then the hypothesis (d1b) is satisfied. 

\medskip

Another example, not of the form (1), is the equation
$$
x'(t)=g(x(t-1-\delta(x(t)+x(t-2))))
$$
from \cite{MVW} with $g:\R\to[-1,1]$ and $\delta:\R\to[-1,1]$ continuously differentiable,
which has the form (2) with  
$n=1$, $r=2$, $K=\{1\}$, $F=\R=W$, $V=\R$,
$$
L\phi=\phi(0)+\phi(-2),\quad\mbox{and}\quad d_1(w)=1+\delta(w).
$$
Obviously the hypothesis (b) holds in this case. A further equation of the form (2) with a single delay but not of the form (1) has been investigated by Kennedy in \cite{Ke2}.

\medskip

Similarly Theorems 3.5, 4.7 and 4.8 cover systems of the form
$$
x'(t)=h(x(t),x(t-\hat{\rho}_1),\ldots,x(t-\hat{\rho}_m))
$$
with a real function $h$ and multiple delays which are given by
$$
\hat{\rho}_j=\rho_j(x(0),x(t-\sigma_j))\quad\mbox{for}\quad j=1,\ldots,m,
$$
with functions $\rho_j:\R^2\to[0,M]$ and reals $\sigma_j\in(0,M]$. These are special cases of the  equations studied by Mallet-Paret, Nussbaum, and Paraskevopoulos
in \cite{M-PNP}. Theorems 4.7 and 4.8 also cover the system studied by Lv, Yuan, Pei, and Li in \cite{LYPL}, for which $U_{\emptyset}=U$.

\medskip

Theorem 4.7 contains the result from \cite[Theorem 5.1]{W6} on solution manifolds which are almost graphs over $X_0$ as a special case. In order to see this we recall that 
in \cite{W6} the system
$$
x_{\nu}'(t)=g^{\ast}_{\nu}(x_1(t),\ldots,x_n(t),x_1(t-\delta_1(Lx_t)),\ldots,x_n(-\delta_n(Lx_t))),\quad\nu\in\{1,\ldots,n\},
$$ 
is considered, for $L:C_n\to F$ linear and continuous, $\dim\,F<\infty$, $W\subset F$ open, and the maps $\delta_{\nu}:W\to(0,r]$, $\nu=1,\ldots,n$, and $g^{\ast}:\R^{2n}\supset V^{\ast}\to\R^n$ continuously differentiable, under the assumption that there exist $\phi\in C^1_n$ with $L\phi\in W$ and
$(\phi(0),\phi_1(-\delta_1(L\phi)),\ldots,\phi_n(-\delta_n(L\phi)))\in V^{\ast}$. It has the form 
$$
x'(t)=f^{\ast}(\phi)
$$
for $f^{\ast}:U^{\ast}\to\R^n$ given by
\begin{eqnarray*}
U^{\ast} & = & \{\phi\in C^1_n:L\phi\in W\,\,\mbox{and}\,\,(\phi(0),\phi_1(-\delta_1(L\phi)),\ldots,\phi_n(-\delta_n(L\phi)))\in V^{\ast}\},\\	
f^{\ast}(\phi) & = & g^{\ast}(\phi(0),\phi_1(-\delta_1(L\phi)),\ldots,\phi_n(-\delta_n(L\phi))),	
\end{eqnarray*}
and is Eq. (2) for ${\mathbf k}=n+1$ and
\begin{eqnarray*}
d_1 & = & 0,\\
d_{\nu} & = & \delta_{\nu-1}\quad\mbox{for}\quad\nu=2,\ldots,n+1,\\
V & = & \{\xi\in(\R^n)^{n+1}:(\xi_1,\xi_{21},\ldots,\xi_{n+1,n})\in V^{\ast}\},\\
g(\xi_1,\ldots,\xi_{n+1}) & = & g^{\ast}(\xi_1,\xi_{21},\ldots,\xi_{n+1,n}),
\end{eqnarray*}
with the components $\xi_{\mu\nu}\in\R$ of $\xi_{\mu}\in\R^n$, for $\mu=1,\ldots,n+1$ and $\nu=1,\ldots,n$. We have
\begin{eqnarray*}
U & = & \{\phi\in C^1_n:L\phi\in W\,\,\mbox{and}\,\,(\phi(-d_1(L\phi)),\phi(-d_2(L\phi)),\ldots,\phi(-d_{n+1}(L\phi)\in V\}\\	
& = & \{\phi\in C^1_n:L\phi\in W\,\,\mbox{and}\,\,(\phi(0),\phi_1(-\delta_1(L\phi)),\ldots,\phi_n(-\delta_n(L\phi))\in V^{\ast}\}\\	
& = & U^{\ast}	
\end{eqnarray*}
and $f(\phi)=f^{\ast}(\phi)$ on $U=U^{\ast}$. The equation $d_1=0$ and positivity of the delay functions $d_{\nu}=\delta_{\nu-1}$ for $2\le\nu\le n+1$ combined yield
$$
U=U_J\quad\mbox{for}\quad J=\{1\}\subset\{1,\ldots,n+1\}=K,
$$
so that the hypothesis (J) in Theorem 4.7 is satisfied.

\medskip

One may ask for conditions which ensure that a solution manifold is a graph over $X_0$. In the sequel this is guaranteed only by Theorem 3.5 in case $U_K=U$, where all delays are zero and Eq. (2) reduces to an ordinary differential equation.
For systems in the general form (4) results from \cite[Section 2]{W6} show that the solution manifold $X_f$ is a graph over $X_0$ provided
$f:C^1_n\supset dom \to\R^n$  is continuously differentiable, condition (e) is satisfied,  and there exists
$z\in(0,r)$ so that
$$
f(\phi)=f(\psi)\,\,\mbox{for all}\,\,\phi\in dom\,\,\mbox{and}\,\,\psi\in dom \,\,\mbox{with}\,\,\phi(t)=\psi(t)\,\,\mbox{on}\,\,[-r,z]\cup\{0\}.
$$
The preceding property generalizes the condition that in a system 
$$
x'(t)=g(x(t), x(t-\Delta_2(x_t)),\ldots,x(t-\Delta_{{\bf k}}(x_t))),
$$
analogous to Eq. (2) but with more general discrete delays $\Delta_k:C^1_n\supset dom\to[0,r]$, $2\le k\le{\bf k}$, all of them are bounded away from zero. The results from \cite[Section 2]{W6} apply to the system studied in \cite{LYPL}, for example, and show that in addition to the almost graph representation obtained via Theorem 4.7 the solution manifold associated with the system from \cite{LYPL} also is a graph over $X_0$.

\medskip

Another possibility to write a solution manifold $X_f$ as a graph, for continuously differentiable $f:C^1_n\supset dom\to\R^n$ with property (e) and $D_ef:dom\to L_c(C_n,\R^n)$ bounded, is due to Krisztin and Rezounenko \cite{K,KR}, see the proof of \cite[Lemma 1]{KR}. 

\medskip

Solution manifolds which are graphs occur also
in \cite{QW,KW}.

\medskip

{\bf Notation, preliminaries.} For subsets $A,B$ of a vectorspace $V$ over the field $\R$ and for $c\in V$ we use the abbreviations
\begin{eqnarray*}
	A+c & = & \{v\in V:\mbox{There exists}\,\,a\in A\,\,\mbox{with}\,\,v=a+c\},\\
	A+B & = & \{v\in V:\mbox{There exist}\,\,a\in A\,\,\mbox{and}\,\,b\in B\,\,\mbox{with}\,\,v=a+b\},\\
	\R\,c & = & \{v\in V:\mbox{There exists}\,\,\xi\in \R\,\,\mbox{with}\,\,v=\xi c\}.
\end{eqnarray*}

\medskip

Finite-dimensional vectorspaces are always equipped with the canonical topology which makes them topological vectorspaces. 

\medskip

An upper index as in $(x_1,\ldots,x_n)^{tr}\in\R^n$ denotes the transpose of the row vector $(x_1,\ldots,x_n)$. Vectors in $\R^n$ which occur as argument of a map are always written as row vectors. The vectors of the canonical basis of $\R^n$ are denoted by $e_{\nu}$, $\nu\in\{1,\ldots,n\}$; $e_{\nu\mu}=1$ for $\nu=\mu$ and $e_{\nu\mu}=0$ for $\nu\neq\mu$, $\nu$ and $\mu$ in $\{1,\ldots,n\}$.   

\medskip

For subsets $A\subset B$ of a topological space $T$ we say $A$ is open in $B$ if $A$ is open with respect to the relative topology on $B$. Analogously for $A$ closed in $B$.

\medskip

The relation $W_1\subset\subset W_2$ for open subsets of a finite-dimensional vectorspace means that the closure $\overline{W_1}$ of $W_1$ is compact and contained in $W_2$.

\medskip

Derivatives  and partial derivatives of a map at a given argument are continuous linear maps, indicated by a capital $D$. In case of real functions on domains in $\R$ and in $\R^n$, $\phi'(t)=D\phi(t)1$ and $\partial_{\nu}g(x)=D_{\nu}g(x)1$, respectively.

\medskip

For $n=1$ we abbreviate $C=C_1$ and $C^1=C^1_1$.

\medskip

Let $C^1_{n\times n}=C^1([-r,0],\R^{n\times n})$
denote the Banach space of continuously differentiable maps from $[-r,0]$ into the vectorspace $\R^{n\times n}$ of $n\times n$-matrices with real entries which is
analogous to the space $C^1_n$.
For a map $A\in C^1_{n\times n}$ its components in $C^1$ are given by $A_{\nu\mu}(t)=(A(t))_{\nu\mu}\in\R$, for $\nu,\mu$ in $\{1,\ldots,n\}$ and $t\in[-r,0]$. The columns of $A$ are the maps $A_{\mu}\in C^1_n$, $\mu\in\{1,\ldots,n\}$, with components 
$A_{1\mu},\ldots, A_{n\mu}$ in $C^1$.

\medskip

We define the product of a function $\phi\in C$ with a vector $q=(q_1,\ldots,q_n)^{tr}\in\R^n$ as the map $\phi\cdot q\in C_n$ with the components $q_1\phi,\ldots,q_n\phi$ in $C$. For  $A\in C^1_{n\times n}$ and $q\in\R^n$ we define their product by
$$
A\cdot q=\sum_{\mu=1}^nq_{\mu}A_{\mu}\in C^1_n,
$$
and
$$
A\cdot\R^n=\{\phi\in C^1_n:\mbox{There exists}\,\, q\in\R^n\,\,\mbox{with}\,\,\phi=A\cdot q\}.
$$

The evaluation map
$$
ev:C\times[-r,0]\ni(\phi,t)\mapsto\phi(t)\in\R
$$
is continuous but not locally Lipschitz. The restricted evaluation map
$$
Ev=ev|(C^1\times(-r,0))
$$ 
is continuously differentiable with
$$
D\,Ev(\phi,t)(\chi,s)=D_1Ev(\phi,t)\chi+D_2Ev(\phi,t)s=\chi(t)+\phi'(t)s
$$
for all $\phi\in C^1$, $t\in(-r,0)$, $\chi\in C^1$, $s\in\R$.

\medskip

Differentiation $\partial:C^1_n\ni\phi\mapsto \phi'\in C_n$ and  
evaluation  
$ev_0:C_n\ni\phi\mapsto\phi(0)\in\R^n$ at $t=0$ are  continuous linear maps.

\medskip

The tangent cone of a subset $X$ of a Banach space $B$ at a point $x\in X$  is the set $T_xX$ of all tangent vectors $c'(0)=Dc(0)1\in B$ of continuously differentiable curves $c:J\to B$, with $J\subset\R$ an open interval, $0\in J$,  
$c(J)\subset X$, and $c(0)=x$.

\medskip

The tangent spaces of a solution manifold $X_f=\{\phi\in U:\phi'(0)=f(\phi)\}\neq\emptyset$ associated with a continuously differentiable map $f:U\to\R^n$ with property (e) are given by
$$
T_{\phi}X_f=\{\chi\in C^1_n:\chi'(0)=Df(\phi)\chi\},\quad \phi\in X_f.
$$

In the sequel a diffeomorphism is a continuously differentiable injective map with open image whose inverse is continuously differentiable.

\medskip

\section{The solution manifold}

In the following we consider $r,n,F,L,{\bf k},K,W,d_k$ for $k\in K,V,$ and $g$ as in Section 1, with property (V). The set $U$ is given by Eq. (3), and the map $f$ is given by Eq. (5).

\begin{prop}
$U$ is an open subset of $C^1_n$, $f$ is continuously differentiable with property (e), and for every $\mu\in\{1,\ldots,n\}$ and $\phi\in U$ and $\chi\in C^1_n$,
\begin{equation}
Df_{\mu}(\phi)\chi= \sum_{k=1}^{{\bf k}}\sum_{\nu=1}^n\partial_{(k-1)n+\nu}g_{\mu}(\widehat{\phi})[	\chi_{\nu}(-d_k(L\phi))-	\phi_{\nu}'(-d_k(L\phi))Dd_k(L\phi)L\chi].
\end{equation}	
\end{prop}

{\bf Proof.}
1. On the set $U\neq\emptyset$. As the restriction of $L$ to $C^1_n$ is continuous we see that the  preimage 
$\{\phi\in C^1_n:L\phi\in W\}=L^{-1}(W)\cap C^1_n$ 
is open in $C^1_n$. Consider the map from this preimage into $(\R^n)^{{\bf k}}$ with the $n{\bf k}$ components given by
$$
\phi_{\nu}(-d_k(L\phi))=ev(\phi_{\nu},-d_k(L\phi)),
$$
for
$\nu\in\{1,\ldots,n\}$ and $k\in K$.
Each of these components is continuous. Consequently the map
$$
L^{-1}(W)\cap C^1_n\ni\phi\mapsto\widehat{\phi}\in(\R^n)^{{\bf k}}
$$
is continuous, and it follows that 
\begin{eqnarray*}
U & = & \{\phi\in C^1_n:L\phi\in W\quad\mbox{and}\quad\widehat{\phi}\in V\}\\
& = & \{\phi\in L^{-1}(W)\cap C^1_n:\widehat{\phi}\in V\}
\end{eqnarray*}
is open (in $L^{-1}(W)\cap C^1_n$, in $C^1_n$).

\medskip

2. Let $\mu\in\{1,\ldots,n\}$. Proof that the component $f_{\mu}:C^1_n\supset U\to\R$ is continuously differentiable.

\medskip

2.1. Let  $\nu\in\{1,\ldots,n\}$ and $k\in K$. We proceed to show that the map $h_{\nu k}:C^1_n\supset U\to\R$ given by $h_{\nu k}(\phi)=\phi_{\nu}(-d_k(L\phi))$ is continuously differentiable. 
In order to circumvent talking about differentiability of the evaluation map $C^1\times[-r,0]\ni(\phi,t)\mapsto \phi(t)\in\R$, which is not defined on an open subset of a Banach space, we consider the continuous linear map
$$
C^1\ni\phi\mapsto\tilde{\phi}\in C^1((-r-1,1),\R)
$$
given by $\tilde{\phi}(t)=\phi(t)$ on $[-r,0]$,
$\tilde{\phi}(t)=\phi(-r)+(t+r)\phi'(-r)$ on $[-r-1,-r]$, and 
$\tilde{\phi}(t)=\phi(0)+t\phi'(0)$ on $[0,1]$, and use that the evaluation map 
$$
Ev:C^1([-r-1,1],\R)\times(-r-1,1)\ni(\psi,t)\mapsto\psi(t)\in\R
$$
(which is defined on an open subset of the space $C^1([-r-1,1],\R)\times\R$) is continuously differentiable with
$$
D\,Ev(\psi,t)(\chi,s)=\chi(t)+\psi'(t)s.
$$
Notice that for every $\phi\in U\subset C^1_n$ we have
$$
h_{\nu k}(\phi)=Ev(\tilde{\phi_{\nu}},-d_k(L\phi)).
$$
The chain rule applies and yields that $h_{\nu k}$ is continuously differentiable with
$$
Dh_{\nu k}(\phi)\chi=\chi_{\nu}(-d_k(L\phi))-\phi_{\nu}'(-d_k(L\phi))Dd_k(L\phi)L\chi
$$
for every $\phi\in U$ and $\chi\in C^1_n$.

\medskip

2.2 For every $\phi\in U$,
$$
f_{\mu}(\phi)=g_{\mu}(h_{11}(\phi),\ldots,h_{n1}(\phi);\ldots;h_{1{\mathbf k}}(\phi),\ldots,h_{n{\mathbf k}}(\phi)).
$$ 
The chain rule yields that $f_{\mu}$ is continuously differentiable with
$$
Df_{\mu}(\phi)\chi= \sum_{k=1}^{{\bf k}}\sum_{\nu=1}^n\partial_{(k-1)n+\nu}g_{\mu}(\widehat{\phi})[
\chi_{\nu}(-d_k(L\phi))-
\phi_{\nu}'(-d_k(L\phi))Dd_k(L\phi)L\chi]
$$
for every $\phi\in U$ and $\chi\in C^1_n$.

\medskip
 
2.3. The right hand side of the preceding equation also defines linear maps $D_ef_{\mu}(\phi):C_n\to\R$, for $\phi\in U$ and $\mu=1,\ldots,n$.

\medskip

3. From Part 2 we infer that the map $f$ is continuously differentiable, and that Eq. (6) holds. Using the continuity of differentiation $\partial:C^1\to C$ and of the evaluation map $ev:C\times[-r,0]\to\R$ we infer that the map
$$
U\times C_n\ni(\phi,\chi)\mapsto D_ef(\phi)\chi\in\R^n
$$
given by the components
$$
(D_ef(\phi)\chi)_{\mu}=D_ef_{\mu}(\phi)\chi
$$
is continuous, which means that $f$ has property (e). $\Box$

\medskip

The next lemma will be instrumental in the proof of $X_f\neq\emptyset$ and in Sections 3 and 4.  

\begin{lem}
(\cite[Proposition 4.1]{W6})
Let a continuous linear functional $\lambda:C\to\R^q$ and $z\in[-r,0)$ be given. There exists $\phi\in C^1$ with
$\phi'(0)=1$, $\lambda\phi=0$, and $\phi(t)=0$ on $[-r,z]\cup\{0\}$.
\end{lem}

\begin{prop}
For every $J\subset K$ with $U_J\neq\emptyset$ we have $X_{fJ}\neq\emptyset$, and $X_f\neq\emptyset$ is a continuously differentiable submanifold of codimendion $n$ in $C^1_n$.
For every $\phi\in X_f$,
\begin{eqnarray*}
T_{\phi}X_f & = & \{\chi\in C^1_n:\mbox{For every}\quad\mu\in\{1,\ldots,n\},\\
& & (\chi_{\mu})'(0)=
\sum_{k=1}^{{\bf k}}\sum_{\nu=1}^n\partial_{(k-1)n+\nu}g_{\mu}(\widehat{\phi})[
	\chi_{\nu}(-d_k(L\phi))\\
& & -
	\phi_{\nu}'(-d_k(L\phi))Dd_k(L\phi)L\chi]\}.
\end{eqnarray*}  
\end{prop}
   
{\bf Proof.} 1. Let $J\subset K$ with $U_J\neq\emptyset$ be given. Proof of $X_{fJ}\neq\emptyset$. Choose $\phi\in U_J$. Then $L\phi\in W_J$ and $\widehat{\phi}\in V$. Set $q=f(\phi)-\phi'(0)$, and $z=\max_{k\in K\setminus J}(-d_k(L\phi))$ if $J\neq K$, and $z=-r/2$ if $J=K$. Then $-r\le z<0$. 

\medskip

1.1. Let $\nu\in\{1,\ldots,n\}$ be given. We show that there exist $\eta_{\nu}\in C^1$ with  $\eta_{\nu}(0)=0$, $\eta_{\nu}'(0)=1$, $L(\eta_{\nu}\cdot e_{\nu})=0$, and $\eta_{\nu}(-d_k(L\phi))=0$ for all $k\in K$. Proof of this: Choose an isomorphism  $\iota:F\to\R^{\dim\,F}$ and apply Lemma 2.2 to the map $\lambda:C\to\R^{\dim\,F}$  given by $\lambda\psi=\iota L(\psi\cdot e_{\nu})$, with $z$ from above. This yields $\eta_{\nu}\in C^1$ with $\eta_{\nu}'(0)=1$, $L(\eta_{\nu}\cdot e_{\nu})=0$,  and $\eta_{\nu}(t)=0$ on $[-r,z]\cup\{0\}$. It follows that $\eta_{\nu}(-d_j(L\phi))=0$ for every $j\in J$ since $d_j(L\phi)=0$ for $j\in J$ (due to $L\phi\in W_J$). In case $J\neq K$ and $k\in K\setminus J$ we get  $\eta_{\nu}(-d_k(L\phi))=0$ from $-d_k(L\phi)\le z$. 
 
\medskip

1.2. Set
$$
\psi=\sum_1^nq_{\nu}\eta_{\nu}\cdot e_{\nu}.
$$
We show $\phi+\psi\in X_{fJ}$. Observe that $L\psi=0$, $\psi'(0)=q$, and $\psi(-d_k(L\phi))=0$ for all $k\in K$. It follows that $(\phi+\psi)'(0)=\phi'(0)+q=f(\phi)$, and it remains to deduce $\phi+\psi\in U_J$ and $f(\phi+\psi)=f(\phi)$. Using $L(\phi+\psi)=L\phi\in W_J$ we get
\begin{eqnarray*}
\widehat{\phi+\psi} & = & (\ldots,(\phi+\psi)(-d_k(L(\phi+\psi))),\ldots)=(\ldots,(\phi+\psi)(-d_k(L\phi)),\ldots)\\
& = & (\ldots,\phi (-d_k(L(\phi))+\psi(-d_k(L\phi)),\ldots)=(\ldots,\phi(-d_k(L(\phi)),\ldots)\\
& = & \widehat{\phi}\in  V,
\end{eqnarray*}
hence $\phi+\psi\in U_{J}$ and
$$
f(\phi+\psi)=g(\widehat{\phi+\psi})=g(\widehat{\phi})=f(\phi).
$$

2. The hypothesis (V) yields $U_J\neq\emptyset$ for some $J\subset K$. Consequently, $\emptyset\neq X_{fJ}\subset X_f$. As $f$ is continuously differentiable with property (e) we have from remarks in Section 1 that $X_f\neq\emptyset$
is a continuously differentiable submanifold of codimension $n$ in $C^1_n$.

\medskip

3.  The assertion about the tangent spaces follows from $T_{\phi}X_f=\{\chi\in C^1_n:\chi'(0)=Df(\phi)\chi\}$ in combination with the formula for the derivatives of the components of $f$ in Proposition 2.1. $\Box$

\section{On the subset of disappearing delays}

 In this section we assume $U_K\neq\emptyset$. Then $X_{fK}$ is the subset of the solution manifold on which $d_k(L\phi)=0$ for all $k\in K$. By Proposition 2.3, $X_{fK}\neq\emptyset$. The first step on the way to a graph representation of a neighbourhood of $X_{fK}$ in $X_f$, and upon that to an associated map which is better than an almost graph diffeomorphism, is a suitable complementary space for $X_0$ in $C^1_n$.  We choose some  $z\in(-r,0)$ and apply Lemma 2.2 as in the proof of Proposition 2.3. This yields  $\psi_{\nu}\in C^1$, $\nu\in\{1,\ldots,n\}$,  with $L(\psi_{\nu}\cdot e_{\nu})=0$, $\psi_{\nu}(0)=0$, and $\psi_{\nu}'(0)=1$. We define $Y_K\in C^1_{n\times n}$ by its columns $Y_{K\nu}=\psi_{\nu}\cdot e_{\nu}\in C^1_n$, $\nu\in\{1,\ldots,n\}$, set
$$
Q_K=Y_K\cdot\R^n,
$$
and introduce the map $R^K:C^1_n\to C^1_n$ given by
$$
R^K\phi=\phi-Y_K\cdot\phi'(0)=\phi-Y_K\cdot(ev_0\partial\phi).
$$
$R^K$ is linear and continuous. 

\begin{prop}
(i) For every $x\in\R^n$,
\begin{eqnarray}
(Y_K\cdot x)(0) & = & 0\in\R^n,\\
(Y_K\cdot x)'(0) & = & x,\\
L(Y_K\cdot x) & = & 0\in F.
\end{eqnarray}

(ii) $\dim\,Q_K=n$, $C^1_n=H\oplus Q_K$, and $R^K$ is a projection along $Q_K$ onto $X_0$.

\medskip

(iii) For every $\phi\in U_K$ and all $x\in\R^n$,
$$
\phi+Y_K\cdot x\in U_K\quad\mbox{and}\quad f(\phi+Y_K\cdot x)=f(\phi).
$$
\end{prop}
{\bf Proof.} 1. On assertion (i). For every $x\in\R^n$,
$$
Y_K\cdot x=\sum_{\nu=1}^nx_{\nu}Y_{K\nu}=\sum_{\nu=1}^nx_{\nu}\psi_{\nu}\cdot e_{\nu}.
$$
Use  $\psi_{\nu}(0)=0$, $\psi_{\nu}'(0)=1$, and $L(\psi_{\nu}\cdot e_{\nu})=0$ for every $\nu\in\{1,\ldots,n\}$.

\medskip

2. On assertion (ii). For every $\nu\in\{1,\ldots,n\}$, $(Y_{K\nu})'(0)=\psi_{\nu}'(0)e_{\nu}=e_{\nu}$. This implies that the columns $Y_{K\nu}\in C^1_n$ are linearly independent, and it follows that $\dim\,Q_K=n$. For every $\phi\in C^1_n$ we have
$$
(R^K\phi)'(0)=\phi'(0)-(Y_K\cdot\phi'(0))'(0)=\phi'(0)-\phi'(0)=0,
$$
hence $R^KC^1_n\subset X_0$. Also, for $\phi\in X_0$, $R^K\phi=\phi$. We infer $R^KC^1_n=X_0$ and $R^K\circ R^K=R^K$. So $R^K$ is a projection onto $X_0$, with
$$
(R^K)^{-1}(0)=\{\phi\in C^1_n:\phi=Y_K\cdot\phi'(0)\}\subset Q_K.
$$
For every $\phi\in Q_K$ we have $\phi=Y_K\cdot x$ for some $x\in\R^n$, and Eq. (8) yields $\phi'(0)=x$, 
hence $\phi=Y_K\cdot\phi'(0)$, or, $\phi\in(R^K)^{-1}(0)$. It follows that $(R^K)^{-1}(0)=Q_K$ and
$C^1_n=R^KC^1_n\oplus (R^K)^{-1}(0)=X_0\oplus Q_K$. 

\medskip

3. On assertion (iii). Let $\phi\in U_K$ and $x\in\R^n$ be given. Using assertion (i) we get
$(\phi+Y_K\cdot x)(0)=\phi(0)$ and $L(\phi+Y_K\cdot x)=L\phi\quad(\in W)$. Hence $L(\phi+Y_K\cdot x)\in W$, and for every $k\in K$, $d_k(L(\phi+Y_K\cdot x))=d_k(L\phi)=0$ (since $\phi\in U_K$). Furthermore, for every $k\in K$,
\begin{eqnarray*}
(\phi+Y_K\cdot x)(-d_k(L(\phi+Y_K\cdot x))) & = & (\phi+Y_K\cdot x)(-d_k(L(\phi))\\
& = & (\phi+Y_K\cdot x)(0)=\phi(0)\quad\mbox{(with Eq. (7))}\\
& = & \phi(-d_k(L\phi)), 
\end{eqnarray*}
which yields $\widehat{\phi+Y_K\cdot x}=\widehat{\phi}\quad(\in V)$. Altogether, we get $\phi+Y_K\cdot x\in U_K$ and
$$
f(\phi+Y_K\cdot x)=g(\widehat{\phi+Y_K\cdot x})=g(\widehat{\phi})=f(\phi).\qquad\Box
$$

\begin{prop}
The restriction of $R^K$ to $X_{fK}$ is injective, and for every $\phi\in X_{fK}$ the map $R^K$ defines an isomorphism from the tangent space $T_{\phi}X_f$ onto $X_0$.
\end{prop}

{\bf Proof.} 1. On injectivity of $R^K|X_{fK}$. Let $\phi$ and $\psi$ in $X_{fK}=X_f\cap U_K$ be given with $R^K\phi=R^K\psi$. Then
$$
\phi-Y_K\cdot f(\phi)=\phi-Y_K\cdot\phi'(0)=\psi-Y_K\cdot\psi'(0)=\psi-Y_K\cdot f(\psi),
$$
and in order to obtain $\phi=\psi$ it remains to prove $f(\phi)=f(\psi)$. We have $\phi\in X_{fK}\subset U_K$ and 
$$
\psi=\phi+(\psi-\phi)=\phi+Y_K\cdot(f(\psi)-f(\phi)).
$$
Part (iii) of Proposition 3.1 yields $f(\psi)=f(\phi)$.

\medskip

2. Let $\phi\in X_{fK}\subset X_f$ and $\chi\in T_{\phi}X_f$ be given.

\medskip

2.1. In order to obtain injectivity of the restriction of $R^K$ to $T_{\phi}X_f$ we consider $\chi\in T_{\phi}X_f$ with $R^K\chi=0$ and show $\chi=0$. From $R^K\chi=0$, $\chi=Y_K\cdot\chi'(0)$. Using Eq. (9) we infer, $L\chi=0$. By $\phi\in U_K$, $d_k(L\phi)=0$ for every $k\in K$. Using this and Eq. (7) we obtain that for every $\nu\in\{1,\ldots,n\}$ and for every $k\in K$,
$$
\chi_{\nu}(-d_k(L\phi))=\chi_{\nu}(0)=(Y_K\cdot\chi'(0))_{\nu}(0)=0.
$$
Now the representation of $T_{\phi}X_f$ from Proposition 2.3 yields $\chi'(0)=0$. Hence $\chi=Y_K\cdot\chi'(0)=0$.

\medskip

2.2. Proof of $X_0\subset R^KT_{\phi}X_f$. Let $\chi\in X_0$ be given. Consider $x\in\R^n$
with the components
$$
x_{\mu}=\sum_{k=1}^{{\bf k}}\sum_{\nu=1}^n\partial_{(k-1)n+\nu}g_{\mu}(\widehat{\phi})[\chi_{\nu}(0)-\phi_{\nu}'(0)Dd_k(L\phi)L\chi],\quad\mu\in\{1,\ldots,n\},
$$
and
$$
\eta=\chi+Y_K\cdot x\in C^1_n.
$$
We have $d_k(L\phi)=0$ for all $k\in K$ due to $\phi\in U_K$, and $(Y_K\cdot x)(0)=0$
due to Eq. (7), and $(Y_K\cdot x)'(0)=x$ due to Eq. (8), and $L(Y_K\cdot x)=0$ due to Eq. (9). Using these equations we obtain
that for every $\mu\in\{1,\ldots,n\}$,
\begin{eqnarray*}
(\eta_{\mu})'(0) & = & (\chi_{\mu})'(0)+x_{\mu}\\
& = & x_{\mu}\quad\mbox{(since}\quad\chi\in X_0,\chi'(0)=0)\\
& = & \sum_{k=1}^{{\bf k}}\sum_{\nu=1}^n\partial_{(k-1)n+\nu}g_{\mu}(\widehat{\phi})[(\chi_{\nu}+(Y_K\cdot x)_{\nu})(-d_k(L\phi))\\
& & -\phi_{\nu}'(-d_k(L\phi))Dd_k(L\phi)L(\chi+Y_K\cdot x)]\\
& = & \sum_{k=1}^{{\bf k}}\sum_{\nu=1}^n\partial_{(k-1)n+\nu}g_{\mu}(\widehat{\phi})[\eta_{\nu}(-d_k(L\phi))\\
& & -\phi_{\nu}'(-d_k(L\phi))Dd_k(L\phi)L\eta],
\end{eqnarray*}
which means $\eta\in T_{\phi}X_f$, according to the representation of $T_{\phi}X_f$ from Proposition 2.3. Also,
\begin{eqnarray*}
R^K\eta & = &\eta-Y_K\cdot\eta'(0)=\chi+Y_K\cdot x-Y_K\cdot (\chi+Y_K\cdot x)'(0)\\
& = & \chi+Y_K\cdot x-Y_K\cdot(0+x)\quad\mbox{(with}\quad\chi'(0)=0\quad\mbox{and Eq. (8))}\\
& = & \chi.
\end{eqnarray*}

\medskip

3. The results from Part 2 in combination with $R^KC^1_n\subset X_0$ yield that $R^K$ defines an isomorphism from $T_{\phi}X_f$ onto $X_0$, for $\phi\in X_{fK}$. $\Box$

\begin{cor}
For every $\phi\in X_{fK}$ there is an open neighbourhood $N_{\phi}$ of $\phi$ in $X_f$ so that $R^K$ defines a diffeomorphism from $N_{\phi}$ onto the open neighbourhood $R^KN_{\phi}$ of $R^K\phi$ in $X_0$.
\end{cor}

{\bf Proof.} Let $\phi\in X_{fK}$ be given.  Proposition 3.2 in combination with the Open Mapping Theorem guarantees that $DR^K(\phi)=R^K$ defines a topological isomorphism from $T_{\phi}X_f$ onto $X_0$.  Therefore the Inverse Mapping Theorem yields an open neighbourhood $N_{\phi}$ of $\phi$ in $X_f$ so that $R^K$ defines a diffeomorphism from $N_{\phi}$ onto the open neighbourhood $R^KN_{\phi}$ of $R^K\phi$ in $X_0$. $\Box$
 
\begin{prop}
Suppose that in addition to $U_K\neq\emptyset$ one of the hypotheses (b), (d1b), (K) is satisfied.

\medskip

(i) Then every $\chi\in R^KX_{fK}$ has an open neighbourhood $V_{\chi}$ in $X_0$ so that the restriction of $R^K$ to the set
$$
X_f\cap (R^K)^{-1}(V_{\chi})
$$
is injective,

\medskip

and

\medskip

(ii) the restriction of $R^K$ to the set
$$
X_f\cap(\cup_{\chi\in R^KX_{fK}}(R^K)^{-1}(V_{\chi}))
$$
is injective.
\end{prop}

{\bf Proof.} 1. Proof of assertion (i).  Suppose the assertion is false. Then there are $\chi=R^K\phi$ with $\phi\in X_{fK}$ and sequences of elements $\psi_m\neq\tilde{\psi}_m$, $m\in\N$, in $X_f\cap U^K$ with $\chi_m=R^K\psi_m=R^K\tilde{\psi}_m\to\chi$ as $m\to\infty$. For every $m\in\N$,
$$
\psi_m=\chi_m+Y_K\cdot(\psi_m)'(0)=\chi_m+Y_K\cdot f(\psi_m)
$$
and analogously
$\tilde{\psi}_m=\chi_m+Y_K\cdot f(\tilde{\psi}_m)$. By Eq. (7),
$\psi_m(0)=\chi_m(0)=\tilde{\psi}_m(0)$ for every $m\in\N$.

\medskip

1.2. We show that there are a strictly increasing sequence $(m_p)_1^{\infty}$ of positive integers and $x,\tilde{x}$ in $\R^n$ so that $\psi_{m_p}\to \chi+Y_K\cdot x$ and $\tilde{\psi}_{m_p}\to \chi+Y_K\cdot\tilde{x}$ as $p\to\infty$.

\medskip

1.2.1. In case (b) holds the map $g$ is bounded, and $f$ is bounded. The Bolzano-Weierstra{\ss} Theorem yields a convergent subsequence $(f(\psi_{\alpha(m)}))_1^{\infty}$, with  $\alpha:\N\to\N$ strictly increasing, and upon that a convergent subsequence of $(f(\tilde{\psi}_{\alpha(m)}))_1^{\infty}$ given by a strictly increasing map $\beta:\N\to\N$. Set
$m_p=\alpha(\beta(p))$ for $p\in\N$ and let $x\in\R^n$ and $\tilde{x}\in\R^n$ denote the limits of the subsequences
$(f(\psi_{m_p}))_1^{\infty}$ and $(f(\tilde{\psi}_{m_p}))_1^{\infty}$, respectively.
By continuity we get that
$$
\psi_{m_p}=\chi_{m_p}+Y_K\cdot f(\psi_{m_p})
\to\chi+Y_K\cdot x
$$
as $p\to\infty$. Analogously,
$\tilde{\psi}_{m_p}\to\chi+Y_K\cdot\tilde{x}$ as $p\to\infty$.

\medskip

1.2.2. In case (d1b) holds we have $d_1(w)=0$ on $W$, and for every bounded set $B\subset\R^n$ the set $\,\,\{g(v)\in\R^n:v\in V\,\,\mbox{and}\,\,(v_1,\ldots,v_n)^{tr}\in B\}\,\,$ is bounded. 
From $\lim_{m\to\infty}\chi_m=\chi$ and $\psi_m(0)=\chi_m(0)$ for all $m\in\N$ we get that the set 
$$
B=\{\psi_m(0)\in\R^n:m\in\N\}
$$
is bounded. Using $d_1(w)=0$ on $W$ we see that for every $m\in\N$
\begin{eqnarray*}
f(\psi_m) & = & g(\widehat{\psi_m})\\	
& = & g(\psi_m(-d_1(L\psi_m)),\psi_m(-d_2(L\psi_m)),\ldots,\psi_m(-d_{{\bf k}}(L\psi_m)))\\
& = & g(\psi_m(0),\psi_m(-d_2(L\psi_m)),\ldots,\psi_m(-d_{{\bf k}}(L\psi_m)))
\end{eqnarray*}
is contained in the set $\{g(v)\in\R^n:v\in V\,\,\mbox{and}\,\,(v_1,\ldots,v_n)^{tr}\in B\}$, which is bounded due to the hypothesis. Analogously the set  
$\{f(\tilde{\psi}_m)\in\R^n:m\in\N\}$ is bounded.
Proceed as in Part 1.2.1, beginning with applications of the Bolzano-Weierstra{\ss} Theorem. 

\medskip

1.2.3. In case $U_K=U$ we get $\psi_m\in U^K\subset U=U_K$ for every $m\in\N$, and Proposition 3.1 (iii) yields $\chi_m=\psi_m-Y_K\cdot\psi_m'(0)\in U_K$ and $f(\chi_m)=f(\psi_m)=\psi_m'(0)$ for every $m\in\N$. Analogously,
$\chi=R^K\phi=\phi-Y_K\cdot\phi'(0)$ with $\phi\in X_{fK}\subset U_K$ belongs to $U_K=U$. By continuity, $f(\chi_m)\to f(\chi)$ as $m\to\infty$. It follows that
$$
\psi_m=\chi_m+Y_K\cdot f(\psi_m)\to\chi+Y_K\cdot f(\chi)
$$
as $m\to\infty$.
Analogously, $\tilde{\psi}_m\to\chi+Y_K\cdot f(\chi)$ as $m\to\infty$.
Set $m_p=p$ for all $p\in\N$ and $x=\tilde{x}=f(\chi)\in\R^n$.

\medskip

1.3. For $\psi=\chi+Y_K\cdot x$ we show $\psi\in X_{fK}$ and $R^K\psi=\chi$. 

\medskip

1.3.1.  Proof of $\psi\in U_K$. By Eqs. (7) and (9), $\psi(0)=\chi(0)=\phi(0)$ and $L\psi=L\chi=L\phi\in W$.
From $\phi\in X_{fK}\subset U_K$ we obtain $L\phi\in W_K$, so for every $k\in K$, $d_k(L\psi)=d_k(L\phi)=0$, hence
$$
\psi(-d_k(L\psi))=\psi(0)=\phi(0)=\phi(-d_k(L\phi)),
$$
or $\widehat{\psi}=\widehat{\phi}\in V$, and
thereby, $\psi\in U$. As $L\psi=L\phi\in W_K$ we also obtain $\psi\in U_K$. 

\medskip	

1.3.2. Using $X_f\ni\psi_{m_p}\to\psi\in U_K$ for $p\to\infty$ and the fact that $X_f$ is closed in $U$ we infer $\psi\in X_f\cap U_K=X_{fK}$. Using Eqs. (9) and (8), and $\chi\in X_0$, we obtain from the definition $\psi=\chi+Y_K\cdot x$ that $ L\psi=L\chi$ and $\psi'(0)=\chi'(0)+x=x$. Hence $R^K\psi=\psi-Y_K\cdot\psi'(0)=\psi-Y_K\cdot x=\chi$.

\medskip

1.4. Set $\tilde{\psi}=\chi+Y_K\cdot\tilde{x}$. Part 1.3 with $\tilde{x}$ in place of $x$ yields $\tilde{\psi}\in X_{fK}$ and $R^K\tilde{\psi}=\chi$. By Proposition 3.2 the restriction of $R^K$ to $X_{fK}$ is injective, hence
$\psi=\tilde{\psi}=\phi$. For $p\in\N$ sufficiently large we obtain $\psi_{m_p}\in N_{\phi}$ and $\tilde{\psi}_{m_p}\in N_{\phi}$, and arrive at a contradiction to the injectivity of $R^K$ on $N_{\phi}$. 

\medskip

2. On assertion (ii). Let $\phi$ and $\psi$ in 
$$
X_f\cap(\cup_{\chi\in R^KX_{fK}}(R^K)^{-1}(V_{\chi}))
$$
be given with $R^K\phi=R^K\psi$. $R^K\phi$ is contained in $V_{\chi}$ for some $\chi\in R^KX_{fK}$, and $R^K\psi=R^K\phi$ is contained in the same set $V_{\chi}$. Or, both $\phi$ and  $\psi$ belong to
$X_f\cap (R^K)^{-1}(V_{\chi})$.  Part (i) yields $\phi=\psi$. $\Box$

\begin{thm}
Suppose $U_K\neq\emptyset$ and that one of the hypotheses (b), (d1b), (K) is satisfied. 

(i) Then there exists an open neighbourhood $N_K$ of $X_{fK}$ in $X_f$ so that the projection $R^K$ defines a diffeomorphism $R_K$ from $N_K$ onto the open subset $R^KN_K$ of $X_0$, and $N_K$ is a graph over $X_0$.

(ii) The continuously differentiable map $A_K:{\mathcal O}_K\to C^1_n$ given by
$$
{\mathcal O}_K=(R^K)^{-1}(R^KN_K)\quad\mbox{and}\quad A_K(\phi)=\phi-(id-R^K)(R_K^{-1}(R^K\phi))
$$
defines a diffeomorphism onto ${\mathcal O}_K$ with
$$
A_K(X_f\cap{\mathcal O}_K)=X_0\cap A_K({\mathcal O}_K)
$$
and
$$
A_K(\phi)=R^K\phi\quad\mbox{on}\quad X_f\cap{\mathcal O}_K=N_K,\quad A_K(\chi)=\chi\quad\mbox{on}\quad X_0\cap( X_f\cap{\mathcal O}_K).
$$
\end{thm}

The map $A_K$ is an almost graph diffeomorphism, of course, with associated manifold chart $R_K$.

\medskip

{\bf Proof} of Theorem 3.5. 1. On assertion (i).  Recall the neighbourhoods $N_{\phi}$ from Corollary 3.3 and the neighbourhoods $V_{\chi}$ from Proposition 3.4.  The set
$$
N_K=\cup_{\phi\in X_{fK}}(N_{\phi}\cap(R^K)^{-1}(V_{R^K\phi}))
$$
is open in $X_f$ and contains $X_{fK}$. Due to Proposition 3.4 (ii)  the restriction of $R^K$ to $N_K$ is injective.
Corollary 3.3 shows that locally this restriction is given by diffeomorphisms onto open subsets of the space $X_0$. It follows easily  that $R^K$ defines a diffeomorphism $R_K$ from the open subset $N_K$ of $X_f$ onto the open subset $R^KN_K$ of $X_0$. 

Using that $R^K$ is a projection onto $X_0$ one finds that $N_K$ is the set of all
$$
\phi=R_K^{-1}(\chi)=R^ K(R_K^{-1}(\chi))+(id-R^K)(R_K^{-1}(\chi)),\quad\chi\in R^KN_K,
$$
with $(id-R^K)(R_K^{-1}(\chi)) \in Q_K$. This means that $N_K$ is a graph over $X_0$. 

\medskip

2.   Proof of $X_f\cap{\mathcal O}_K=N_K$.  The inclusion $N_K\subset X_f\cap{\mathcal O}_K$ is obvious from $N_K\subset X_f$ and 
$$
N_K\subset(R^K)^{-1}(R^KN_K)={\mathcal O}_K.
$$ 
In order to show the reverse inclusion $X_f\cap{\mathcal O}_K\subset N_K$ let $\phi\in X_f\cap{\mathcal O}_K$ be given. Then 
$$
R^K\phi\in R^KN_K,
$$
or, $R^K\phi=R^K\psi$ for some
$$
\psi\in N_K\subset X_f\cap(\cup_{\chi\in R^KX_{fK}}(R^K)^{-1}(V_{\chi})).
$$
Hence   
$$
R^K\phi=R^K\psi\in V_{\eta}\quad\mbox{for some}\quad\eta\in R^KX_{fK},
$$
or, $\phi\in(R^K)^{-1}(V_{\eta})$.
It follows that both $\phi$ and $\psi$ belong to the set $X_f\cap(R^K)^{-1}(V_{\chi}))$, and $R^K\phi=R^K\psi$.
Using Proposition 3.4 (i) we infer $\phi=\psi\in N_K$.

\medskip

3. The set ${\mathcal O}_K$ is open in $C^1_n$ and the map $A_K:{\mathcal O}_K\to C^1_n$ is continuously differentiable, as well as the map $B_K:{\mathcal O}_K\to C^1_n$ given by $B_K(\psi)=\psi+(id-R^K)(R_K^{-1}(R^K\psi))$. 

\medskip

Proof of $A_K({\mathcal O}_K)\subset {\mathcal O}_K$. Let $\psi=A_K(\phi)$ with $\phi\in{\mathcal O}_K$ be given. Then 
\begin{eqnarray*}
R^K\psi & = & R^K(\phi-(id-R^K)(R_K^{-1}(R^K\phi)))=R^K\phi\\
& & \mbox{(since}\quad R^K\circ(id-R^K)=0)\\
& \in & R^K N_K
\end{eqnarray*} 
which means $A_K(\phi)=\psi\in{\mathcal O}_K$.

\medskip

Analogously, $B_K({\mathcal O}_K)\subset {\mathcal O}_K$. 

\medskip

Proof of $B_K(A_K(\phi))=\phi$ on ${\mathcal O}_K$. For $\phi\in{\mathcal O}_K$ and $\psi=A_K(\phi)$, 
$$
R^K\psi=R^K(\phi-(id-R^K)((R_K)^{-1}(R^K\phi)))=R^K\phi\quad\mbox{(since}\quad R^K\circ(id-R^K)=0),
$$
hence
\begin{eqnarray*}
B_K(A_K(\phi)) & = & \psi+(id-R^K)(R_K^{-1}(R^K\psi))\\
& = & [\phi-(id-R^K)(R_K^{-1}(R^K\phi))]+(id-R^K)(R_K^{-1}(R^K\psi))\\
& = & [\phi-(id-R^K)(R_K^{-1}(R^K\phi))]+(id-R^K)(R_K^{-1}(R^K\phi))=\phi.
\end{eqnarray*}

Analogously, $A_K(B_K(\eta))=\eta$ on ${\mathcal O}_K$. It follows that $A_K$ defines a diffeomorphism from ${\mathcal O}_K$ onto  ${\mathcal O}_K$.

\medskip

4. On $X_f\cap{\mathcal O}_K=N_K$ we have $\phi=R_K^{-1}(R^K\phi)$, hence
$$
A_K(\phi)=\phi-(id-R^K)(R_K^{-1}(R^K\phi))=\phi-(id-R^K)\phi=R^K\phi.
$$
For $\zeta\in X_0\cap(X_f\cap{\mathcal O}_K)$ we infer
$$
A_K\zeta=R^K\zeta=\zeta
$$
since $R^K$ is a projection onto $X_0$.

\medskip

5. Proof of $A_K(X_f\cap{\mathcal O}_K)=X_0\cap A_K({\mathcal O}_K)$. We have $A_K(X_f\cap{\mathcal O}_K)=A_K(N_K)=R^KN_K\subset X_0$ and  $A_K(X_f\cap{\mathcal O}_K)\subset A_K({\mathcal O}_K)$, hence
$A_K(X_f\cap{\mathcal O}_K)\subset X_0\cap A_K({\mathcal O}_K)$. Conversely, for $\chi\in X_0\cap A_K({\mathcal O}_K)$ we have $\chi=R^K\chi$ with $\chi\in A_K({\mathcal O}_K)=
{\mathcal O}_K=(R^K)^{-1}(R^KN_K)$, hence $R^K\chi=R^K\eta$ for some $\eta\in N_K$, and it follows that
$\chi=R^K\chi=R^K\eta=A_K(\eta)\in A_K(N_K)$. Consequently, 
$X_0\cap A_K({\mathcal O}_K)\subset A_K(N_K)=A_K(X_f\cap{\mathcal O}_K)$. $\Box$

\section{Almost graph diffeomorphisms} 

\medskip 

In this section we consider a subset $J\neq K$ of $K$ with $U_J\neq\emptyset$. Then $X_{fJ}\neq\emptyset$, due to Proposition 2.3. We construct a map $R^J:U^J\to C^1_n$ which defines a manifold chart for $X_f$ on an almost graph containing $X_{fJ}$, and obtain finally an almost graph diffeomorphism on a neighbourhood of $X_{fJ}$ in $C^1_n$ which coincides with $R^J$ on the part of $X_f$ in this neighbourhood. In contrast to $R^K$ from Section 3 the map $R^J$ will in general not be a (restricted) projection.

\medskip

The construction of $R^J$ relies on spaces $Q_{J\chi}\subset C^1_n$ complementary to $X_0$, for $\chi\in X_0$ with $L\chi\in W^J$. The next proposition is an adaptation of \cite[Proposition 4.3]{W6} which prepares the  choice of the spaces $Q_{J\chi}$. Its proof makes use of Lemma 2.2. 

\begin{prop}
Let a finite-dimensional vectorspace ${\mathcal F}$, a continuous linear map $\lambda:C\to\R^q$, a non-empty open subset ${\mathcal W}\subset{\mathcal F}$, and a continuous function $d:{\mathcal W}\to(0,r]$ be given. Then there exists a continuously differentiable map $y:{\mathcal W}\to C^1$ such that for every $w\in{\mathcal W}$ we have
\begin{eqnarray}
y(w)(0) & = 0,\\
(y(w))'(0) & = & 1,\\
\lambda(y(w)) & = & 0,\\
& \mbox{and} & \nonumber\\
y(w)(t) & = & 0\quad\mbox{for}\quad-r\le t\le-d(w).
\end{eqnarray}
\end{prop}

{\bf Proof.} 1. Let ${\mathcal F}$, $\lambda:C\to\R^q$, ${\mathcal W}\subset{\mathcal F}$, and $d:{\mathcal W}\to(0,r]$ be given.

\medskip

 2. The case $0<\dim\,{\mathcal F}$.  There are sequences of open subsets $W_{m0},W_{m1},W_m$ of ${\mathcal W}$ with
$$
W_m\subset\subset W_{m+1,0}\quad\mbox{and}\quad \emptyset\neq W_{m0}\subset\subset W_{m1}\subset\subset W_m\quad\mbox{for every}\quad m\in\N
$$
 and ${\mathcal W}=\bigcup_{m\in\N}W_m$. The minima $d_m=\min_{w\in\overline{W_m}}d(w)$ form a decreasing sequence of positive real numbers. For every $m\in\N$ we 
apply Lemma 2.2 (with $z=-d_m/2$) and obtain $\psi_m\in C^1$  with $\lambda\psi_m=0$, $\psi_m(t)=0$ on $[-r,-d_m]\cup\{0\}$, and $\psi_m'(0)=1$. We also choose continuously differentiable functions
$$
a_m:{\mathcal W}\to[0,1]
$$
with $a_m(w)=0$ on $\overline{W_{m0}}$ and $a_m(t)=1$ on ${\mathcal W}\setminus W_{m1}$ for every $m\in\N$, and define
$$
\epsilon_1:W_1\to C^1
$$
by $\epsilon_1(w)=\psi_1+a_1(w)[\psi_2-\psi_1]$, and for every integer $m\ge2$,
$$
\epsilon_m:W_m\setminus\overline{W_{m-1,1}}\to C^1
$$
by $\epsilon_m(w)=\psi_m+a_m(w)[\psi_{m+1}-\psi_m]$. Each map $\epsilon_m$ is continuously differentiable, the domains of $\epsilon_m$ and $\epsilon_p$ with $|m-p|>1$ are disjoint, and for every integer $m\ge2$ we have
$$
\epsilon_m(w)=\psi_m=\epsilon_{m-1}(w)
$$
on the intersection $W_{m-1}\setminus\overline{W_{m-1,1}}$ of the domains of $\epsilon_m$ and $\epsilon_{m-1}$.
The equations
$$
y(w)=\epsilon_1(w)\quad\mbox{on}\quad W_1\quad\mbox{and}\quad y(w)=\epsilon_m(w)\quad\mbox{on}\quad W_m\setminus\overline{W_{m-1,1}}
$$
for integers $m\ge2$ define a continuously differentiable map $y:{\mathcal W}\to C^1$ which obviously satisfies the equations (10)-(12) for every $w\in{\mathcal W}$.

\medskip

 Proof of Eq. (13). For $w\in W_1$, $d(w)\ge d_1\ge d_2$, and we obtain for $-r\le t\le-d(w)$ that $t\le-d_1\le-d_2$. Consequently, $\psi_1(t)=0=\psi_2(t)$, and thereby, 
$$y(w)(t)=\epsilon_1(w)(t)=\psi_1(t)+a_1(w)[\psi_2(t)-\psi_1(t)]=0.
$$ 
For $w\in{\mathcal W}\setminus W_1$ there is an integer $m\ge2$ with
$$
w\in W_m\setminus W_{m-1}\subset W_m\setminus\overline{W_{m-1,1}}.
$$ 
For $t\in[-r,-d(w)]$ we have $-r\le t\le-d(w)\le-d_m\le-d_{m+1}$, hence $\psi_m(t)=0=\psi_{m+1}(t)$, and it follows that
$$
y(w)(t)=\epsilon_m(w)(t)=\psi_m(t)+a_m(w)[\psi_{m+1}(t)-\psi_m(t)]=0.
$$

\medskip

3. The case ${\mathcal F}=\{0\}$. Then ${\mathcal W}=\{0\}$. An application of  Lemma 2.2 with $z=-d(0)$ yields $\phi\in C^1$ with
$\lambda\phi=0$, $\phi'(0)=1$, and $\phi(t)=0$ on $[-r,-d(0)]\cup\{0\}$. Set $y(0)=\phi$. $\Box$

\medskip

In order to apply Proposition 4.1 set ${\mathcal F}=F$ and observe that $W^J\neq\emptyset$, due to $U^J\supset U_J\neq\emptyset$. Set ${\mathcal W}=W^J$ and $d=d^J$ with
$$
d^J:W^J\ni w\mapsto \min_{k\in K\setminus J}d_k(w)\in(0,r].
$$
Let $\nu\in\{1,\ldots,n\}$. Choose an isomorphism $\iota:F\to\R^{\dim\,F}$,  consider the continuous linear map 
$$
\lambda_{\nu}:C\ni\phi\mapsto\iota L(\phi\cdot e_{\nu})\in\R
$$ 
and set $\lambda=\lambda_{\nu}$. Proposition 4.1 yields a continuously differentiable map
$$
y_{J\nu}:W^J\to C^1
$$
with
$$
y_{J\nu}(w)(0)=0,\quad (y_{J\nu}(w))'(0)=1,\quad L(y_{J\nu}(w)\cdot e_{\nu})=0
$$
and
$$
y_{J\nu}(w)(t)=0\quad\mbox{on}\quad[-r,-d^J(w)]
$$
for all $w\in W^J$. The relations
$$
Y_J(w)=(y_{J1}(w)\cdot e_1,\ldots,y_{Jn}(w)\cdot e_n)\in C^1_{n\times n},\quad w\in W^J,
$$
define a continuously differentiable map
$$
Y_J:W^J\to C^1_{n\times n}
$$
with diagonal coefficients of $Y_J(w)$, $w\in W^J$, given by $y_{J\nu}(w)$, $\nu\in\{1,\ldots,n\}$, and
all off-diagonal coefficients equal to $0\in C^1$.

\begin{prop}
(i) For every $w\in W^J$ and $x\in\R^n$,
\begin{eqnarray}
(Y_J(w)\cdot x)(0) & = & 0\in \R^n,\\
(Y_J(w)\cdot x)'(0) & = & x\in \R^n,\\
L(Y_J(w)\cdot x) & = & 0\in F,\\
 \mbox{and for every}\quad k\in K\setminus J, & & \nonumber\\
(Y_J(w)\cdot x)(t) & = & 0\quad\mbox{on}\quad[-r,-d_k(w)]. 
\end{eqnarray}

\medskip

(ii) For all $w\in W^J$, $\tilde{w}\in F$, and $x\in\R^n$,
\begin{eqnarray}
((DY_J(w)\tilde{w})\cdot x)(0) & = & 0\in \R^n,\\
((DY_J(w)\tilde{w})\cdot x)'(0) & = & 0\in \R^n,\\
L((DY_J(w)\tilde{w})\cdot x) & = & 0\in F.
\end{eqnarray}

\medskip

(iii) For every $w\in W^J$, $\dim\,Y_J(w)\cdot\R^n=n$ and 
$$
C^1_n=X_0\oplus Y_J(w)\cdot\R^n.
$$

(iv) For every $\phi\in U_J$ and $x\in\R^n$, 
$$
\phi+Y_J(L\phi)\cdot x\in U_J\quad\mbox{and}\quad f(\phi+Y_J(L\phi)\cdot x)=f(\phi).
$$
\end{prop}

{\bf Proof.} 1. On assertion (i). Let $w\in W^J$ and $x\in\R^n$ be given. The map $Y_J(w)\cdot x\in C^1_n$ has the components $x_{\nu}y_{J\nu}(w)\in C^1$, $\nu\in\{1,\ldots,n\}$. Now the equations (14) and (15) are obvious. For Eq. (16), observe
$$
L(Y_J(w)\cdot x)=L\left(\sum_{\nu=1}^nx_{\nu}y_{J\nu}(w)\cdot e_{\nu}\right)=\sum_{\nu=1}^nx_{\nu}L(y_{J\nu}(w)\cdot e_{\nu})=0.   
$$
Finally, let also $k\in K\setminus J$ and $t\in[-r,-d_k(w)]$ be given. For the components $(x_{\nu}y_{J\nu}(w))(t)$, $\nu\in\{1,\ldots,n\}$, of $(Y_J(w)\cdot x)(t)\in\R^n$ we use $y_{J\nu}(w)(t)=0$ on $[-r,-d^J(w)]\supset[-r,-d_k(w)]$, and obtain Eq. (17).

\medskip

2. On assertion (ii). Let $x\in\R^n$ be given. Due to Eqs. (14)-(16) we have
\begin{eqnarray*}
ev_0(Y_J(w)\cdot x) & = & 0\in\R^n,\\
ev_0\partial(Y_J(w)\cdot x) & = & x\in\R^n,\\
L(Y_J(w)\cdot x) & = & 0\in F
\end{eqnarray*}
for every $w\in W^J$. Use that $ev_0$, $\partial$, and the inclusion map $C^1_n\to C_n$ are linear and continuous, and that the map $Y_J$ is continuously differentiable, and differentiate with respect to $w\in W^J$.

\medskip

3. On assertion (iii). Let $w\in W^J$ be given. The columns $y_{J\nu}(w)\cdot e_{\nu}\in C^1_n$ of $Y_J(w)$ are linearly independent. This yields $\dim\,Y_J(w)\cdot\R^n=n$. For $\phi\in X_0\cap Y_J(w)\cdot\R^n$ we have $\phi'(0)=0$ and
$\phi=\sum_{\nu=1}^nx_{\nu}y_{J\nu}(w)\cdot e_{\nu}=Y_J(w)\cdot x$
with some $x\in\R^n$. Eq. (15) yields $\phi'(0)=x$. Hence $x=0$, and thereby, $\phi=0$. Therefore the sum $X_0+Y_J(w)\cdot\R^n$ is direct. As $X_0$ has codimension $n$ we obtain the desired direct sum decomposition of $C^1_n$.

\medskip

4. On assertion (iv). Let $\phi\in U_J$ and $x\in\R^n$. Then $L\phi\in W_J\subset W^J$. Set $\psi=\phi+Y_J(L\phi)\cdot x$. From Eqs. (14) and (16) we get $\psi(0)=\phi(0)$ and $L\psi=L\phi\quad(\in W_J)$. For $j\in J$ we infer $d_j(L\psi)=d_j(L\phi)=0$, and for $k\in K\setminus J$ we get $d_k(L\psi)=d_k(L\phi)>0$. Consequently,
$$
\psi(-d_j(L\psi))=\psi(0)=\phi(0)=\phi(-d_j(L\phi))\quad\mbox{for every}\quad j\in J.
$$
From Eq. (17) we have $(Y_J(L\phi)\cdot x)(-d_k(L\phi))=0$ for $k\in K\setminus J$ and obtain  
$$
\psi(-d_k(L\psi))=\psi(-d_k(L\phi))=\phi(-d_k(L\phi))\quad\mbox{for}\quad k\in K\setminus J.
$$
It follows that $\widehat{\psi}=\widehat{\phi}\in V$. Hence $\psi\in U_J$ and
$$
f(\psi)=g(\widehat{\psi})=g(\widehat{\phi})=f(\phi).\qquad\Box
$$ 

For each $\chi\in X_0$ 
with $L\chi\in W^J$ the space
$$
Q_{J\chi}=Y_J(L\chi)\cdot\R^n
$$
is complementary to $X_0$ in $C^1_n$, by Proposition 4.2 (iii). Notice that  $LU^J\subset W^J$. The map
$$
R^J:U^J\ni\phi\mapsto\phi-Y_J(L\phi)\cdot\phi'(0)\in C^1_n
$$
has range in $X_0$, due to Eq. (15), and satisfies
$$
R^J(\chi)=\chi\quad\mbox{for every}\quad\chi\in X_0\cap U^J.
$$

\begin{prop}
The map $R^J$ is continuously differentiable with
$$
DR^J(\phi)\chi=\chi-DY_J(L\phi)L\chi\cdot\phi'(0)-Y_J(L\phi)\cdot\chi'(0)
$$
for every $\phi\in U^J$ and $\chi\in C^1_n$.
\end{prop}

{\bf Proof.} The map
$$
B:C^1_{n\times n}\times\R^n\ni(A,x)\mapsto A\cdot x\in C^1_n
$$
is bilinear and continuous, hence continuously differentiable with
$$
DB(A,x)(\tilde{A},\tilde{x})=B(\tilde{A},x)+B(A,\tilde{x})
$$
for all $A\in C^1_{n\times n},x\in\R^n,\tilde{A}\in C^1_{n\times n},\tilde{x}\in\R^n$. For every $\phi\in U^J$ the term $\phi-R^J(\phi)$ equals $B((Y_J\circ L)(\phi),ev_0\partial\phi)$. By the chain rule we infer that $R^J$ is continuously differentiable with
\begin{eqnarray*}
DR^J(\phi)\chi & = & \chi-DB(Y_J(L\phi),\phi'(0))(DY_J(L\phi)L\chi,\chi'(0))\\
& = & \chi-DY_J(L\phi)L\chi\cdot \phi'(0)-Y_J(L\phi)\cdot\chi'(0)
\end{eqnarray*}
for all $\phi\in U^J$ and $\chi\in C^1_n$. $\Box$

\begin{prop}
The restriction of $R^J$ to $X_{fJ}$ 
is injective, and for every $\phi\in X_{fJ}$ the derivative $DR^J(\phi)$ defines an isomorphism from $T_{\phi}X_f$ onto $X_0$.
\end{prop}

{\bf Proof.} 1. On the statement concerning $R^J$. Let $\psi$ and $\phi$ in $X_{fJ}$ be given with $R^J(\psi)=R^J(\phi)$. Then 
$$
\psi-Y_J(L\psi)\cdot f(\psi)=\psi-Y_J(L\psi)\cdot\psi'(0)=\phi-Y_J(L\phi)\cdot\phi'(0)=\phi-Y_J(L\phi)\cdot f(\phi),
$$
and Eq. (16) yields $L\psi=L\phi$. It follows that
$$
\psi=\phi-Y_J(L\phi)\cdot f(\phi)+Y_J(L\psi)\cdot f(\psi)=\phi+Y_J(L\phi)\cdot[f(\psi)-f(\phi)].
$$
We observe $\phi\in X_{fJ}\subset U_J$, apply Proposition 4.2 (iv), and obtain 
$f(\psi)=f(\phi)$, which yields
$\psi=\phi+0=\phi$. 

\medskip

2. On the derivative. Let $\phi\in X_{fJ}$ be given. By Proposition 2.3,
$\chi\in C^1_n$ belongs to $T_{\phi}X_f$ if and only if for every $\mu\in\{1,\ldots,n\}$,
\begin{eqnarray}
(\chi_{\mu})'(0) & = &  \sum_{k=1}^{{\bf k}}\sum_{\nu=1}^n\partial_{(k-1)n+\nu}g_{\mu}(\widehat{\phi})[
\chi_{\nu}(-d_k(L\phi))\nonumber\\
& & -
\phi_{\nu}'(-d_k(L\phi))Dd_k(L\phi)L\chi]\}.
\end{eqnarray}

\medskip

2.1. In order to show injectivity of $DR^J(\phi)$ on $T_{\phi}X_f$ let $\chi\in T_{\phi}X_f$ with $DR^J(\phi)\chi=0$ be given. Using Proposition 4.3 we infer
$$
\chi=DY_J(L\phi)L\chi\cdot\phi'(0)+Y_J(L\phi)\cdot\chi'(0).
$$ 
Using Eqs. (18),(14),(20),(16) we obtain $\chi(0)=0$ and $L\chi=0$. By $L\chi=0$,
$$
\chi=0+Y_J(L\phi)\cdot\chi'(0).
$$
As $\phi\in X_{fJ}\subset U_J$ we have $L\phi\in W_J$, hence $d_j(L\phi)=0$ for $j\in J$, and thereby,
$$
\chi(-d_j(L\phi))=\chi(0)=0\quad\mbox{for all}\quad j\in J.
$$
For $k\in K\setminus J$ we apply Eq. (17) with $L\phi\in W_J\subset W^J$ and get  
$$
\chi(-d_k(L\phi))=(Y_J(L\phi)\cdot\chi'(0))(-d_k(L\phi))=0.
$$
Using $\chi(-d_k(L\phi))=0$ for all $k\in K$ and $L\chi=0$ 
we obtain from Eq. (21) that for every $\mu\in\{1,\ldots,n\}$,
\begin{eqnarray*}
(\chi_{\mu})'(0) & = &  \sum_{k=1}^{{\bf k}}\sum_{\nu=1}^n\partial_{(k-1)n+\nu}g_{\mu}(\widehat{\phi})[
\chi_{\nu}(-d_k(L\phi))\nonumber\\
& & -
\phi_{\nu}'(-d_k(L\phi))Dd_k(L\phi)L\chi]\}=0,	
\end{eqnarray*}
hence $\chi=Y_J(L\phi)\cdot\chi'(0)=0$.

\medskip

2.2. Proof of $DR^J(\phi)T_{\phi}X_f=X_0$. Let $\eta\in X_0$ be given and define $x\in\R^n$ by
\begin{eqnarray}
x_{\mu} & = & \sum_{j\in J}\sum_{\nu=1}^n\partial_{(j-1)n+\nu} g_{\mu}(\widehat{\phi})\{\eta_{\nu}(0)-(\phi_{\nu})'(0)Dd_j(L\phi)L\eta\}\nonumber\\
& & +\sum_{k\in K\setminus J}\sum_{\nu=1}^n\partial_{(j-1)n+\nu} g_{\mu}(\widehat{\phi})\{[\eta+[(DY_J(L\phi)L\eta)\cdot\phi'(0)]_{\nu}(-d_k(L\phi))\nonumber\\
& & -(\phi_{\nu})'(-d_k(L\phi))Dd_k(L\phi)L\eta\}
\end{eqnarray}
for $\mu\in\{1,\ldots,n\}$ and set
\begin{equation}
\chi=\eta+(DY_J(L\phi)L\eta)\cdot\phi'(0)+Y_J(L\phi)\cdot x.
\end{equation}
Then
\begin{eqnarray}
\chi(0) & = & \eta(0)\quad\mbox{(by Eqs. (18) and (14))},\\
\chi'(0) & = & x\quad\mbox{(by Eqs. (19) and (15) and}\quad\eta'(0)=0),\\
L\chi & = & L\eta\quad\mbox{(by Eqs. (20) and (16)).}
\end{eqnarray}
It follows that
\begin{eqnarray*}
DR^J(\phi)\chi & = & \chi-(DY_J(L\phi)L\chi))\cdot\phi'(0)-Y_J(L\phi)\cdot\chi'(0)\\
& = & \chi-(DY_J(L\phi)L\eta))\cdot\phi'(0)-Y_J(L\phi)\cdot
x\\
& = & \eta.	
\end{eqnarray*}
In order to show $\chi\in T_{\phi}X_f$ we verify Eq. (21). Let $\mu\in\{1,\ldots,n\}$ be given. Observe that
\begin{equation}
0=d_j(L\phi)\quad\mbox{for}\quad j\in J,
\end{equation}
due to $\phi\in X_{fJ}\subset U_J$. From Eq. (17) we get
$$
0=(Y_J(L\phi)\cdot x)(-d_k(L\phi))\quad\mbox{for}\quad k\in K\setminus J
$$
which in combination with Eq. (23) yields
\begin{equation}
[\eta_{\nu}+[(DY_J(L\phi)L\eta)\cdot\phi'(0)]_{\nu}](-d_k(L\phi))=\chi_{\nu}(-d_k(L\phi))
\end{equation}
for all $k\in K\setminus J$ and $\nu\in\{1,\ldots,n\}$. We infer
\begin{eqnarray*}
(\chi_{\mu})'(0) & = & x_{\mu}\quad\mbox{(with Eq. (25))}\\
& = & \sum_{j\in J}\sum_{\nu=1}^n\partial_{(j-1)n+\nu} g_{\mu}(\widehat{\phi})\{\chi_{\nu}(-d_j(L\phi))\\
& & -(\phi_{\nu})'(-d_j(L\phi))Dd_j(L\phi)L\chi\}\\
& & \mbox{(with Eqs. (22), (24), (26), (27))}\\
& & +\sum_{k\in K\setminus J}\sum_{\nu=1}^n\partial_{(j-1)n+\nu} g_{\mu}(\widehat{\phi})\{\chi_{\nu}(-d_k(L\phi))\nonumber\\
& & -(\phi_{\nu})'(-d_k(L\phi))Dd_k(L\phi)L\chi\}\\	
& & \mbox{(with Eqs. (22), (28, (26))}\\	
& = & \sum_{k\in K}\sum_{\nu=1}^n\partial_{(k-1)n+\nu} g_{\mu}(\widehat{\phi})\{\chi_{\nu}(-d_k(L\phi))\\
& & -(\phi_{\nu})'(-d_k(L\phi))Dd_k(L\phi)L\chi\},
\end{eqnarray*}
which is Eq. (21). $\Box$

\medskip

As in the proof of Corollary 3.3 we infer the following.

\begin{cor}
For every $\phi\in X_{fJ}$ there is an open neighbourhood $N_{\phi}\subset U^J\cap X_f$ of $\phi$ in $X_f$ so that $R^J$ defines a diffeomorphism from $N_{\phi}$ onto the open neighbourhood $R^J(N_{\phi})$ of $R^J(\phi)$ in $X_0$.	
\end{cor}

\begin{prop}
Suppose that in addition to $U_J\neq\emptyset$ for $J\subset K$ with $J\neq K$ one of the hypotheses (b), (d1b), (J) is satisfied.

\medskip

(i) Then every $\chi\in R^J(X_{fJ})$ has an open neighbourhood $V_{\chi}$ in $X_0$ so that the restriction of $R^J$ to the set
$$
X_f\cap (R^J)^{-1}(V_{\chi})
$$
is injective,

\medskip

and

\medskip

(ii) the restriction of $R^J$ to the subset
$$
X_f\cap(\cup_{\chi\in R^J(X_{fJ})}(R^J)^{-1}(V_{\chi}))
$$
of $X_f\cap U^J$ is injective.
\end{prop}

{\bf Proof.} 1. Suppose assertion (i) is false. Then there are $\chi=R^J(\phi)$ with $\phi\in X_{fJ}$ and sequences of elements $\psi_m\neq\tilde{\psi}_m$, $m\in\N$, in $X_f\cap U^J$ with $\chi_m=R^J(\psi_m)=R^J(\tilde{\psi}_m)\to\chi$ as $m\to\infty$. For every $m\in\N$,
$$
\psi_m=\chi_m+Y_J(L\psi_m)\cdot(\psi_m)'(0)=\chi_m+Y_J(L\psi_m)\cdot f(\psi_m)
$$
and analogously
$\tilde{\psi}_m=\chi_m+Y_J(L\tilde{\psi}_m)\cdot f(\tilde{\psi}_m)$.
By Eq. (14), $\psi_m(0)=\chi_m(0)=\tilde{\psi}_m(0)$ for every $m\in\N$.
By Eq. (16), $L\psi_m=L\chi_m=L\tilde{\psi}_m$ for every $m\in\N$. Using this and $\lim_{m\to\infty}\chi_m=\chi$ and continuity of $L$ we obtain
$$
L\psi_m=L\chi_m\to L\chi\quad\mbox{and}\quad L\tilde{\psi}_m=L\chi_m\to L\chi\quad\mbox{as}\quad m\to\infty.
$$

1.1. We show that there are a strictly increasing sequence $(m_p)_1^{\infty}$ of positive integers and $x,\tilde{x}$ in $\R^n$ so that $\psi_{m_p}\to \chi+Y_J(L\chi)\cdot x$ and $\tilde{\psi}_{m_p}\to \chi+Y_J(L\chi)\cdot\tilde{x}$ as $p\to\infty$.

\medskip

1.1.1. If (b) holds then $g$ is bounded, and $f$ is bounded. The Bolzano-Weierstra{\ss} Theorem yields a convergent subsequence $(f(\psi_{\alpha(m)}))_1^{\infty}$, with  $\alpha:\N\to\N$ strictly increasing, and upon that a convergent subsequence of $(f(\tilde{\psi}_{\alpha(m)}))_1^{\infty}$ given by a strictly increasing map $\beta:\N\to\N$. Set
$m_p=\alpha(\beta(p))$ for $p\in\N$ and let $x\in\R^n$ and $\tilde{x}\in\R^n$ denote the limits of the subsequences
$(f(\psi_{m_p}))_1^{\infty}$ and $(f(\tilde{\psi}_{m_p}))_1^{\infty}$, respectively.
By continuity we get that
$$
\psi_{m_p}=\chi_{m_p}+Y_J(L\psi_{m_p})\cdot f(\psi_{m_p})=
\chi_{m_p}+Y_J(L\chi_{m_p})\cdot f(\psi_{m_p})
\to\chi+Y_J(L\chi)\cdot x
$$
as $p\to\infty$. Analogously,
$\tilde{\psi}_{m_p}\to\chi+Y_J(L\chi)\cdot\tilde{x}$ as $p\to\infty$.

\medskip

1.1.2. If (d1b) holds then $d_1(w)=0$ on $W$ and
for every bounded set $B\subset\R^n$ the set $\,\,
\{g(v)\in\R^n:v\in V\,\,\mbox{and}\,\,(v_1,\ldots,v_n)^{tr}\in B\}\,\,$
is bounded. From $\lim_{m\to\infty}\chi_m=\chi$ and $\psi_m(0)=\chi_m(0)$ for all $m\in\N$ we get that the set
$$
B_0=\{\psi_m(0)\in\R^n:m\in\N\}
$$
is bounded. Using $d_1(w)=0$ on $W$ we see that for every $m\in\N$
\begin{eqnarray*}
f(\psi_m) & = & g(\widehat{\psi_m})\\	
& = & g(\psi_m(-d_1(L\psi_m)),\psi_m(-d_2(L\psi_m)),\ldots,\psi_m(-d_{{\bf k}}(L\psi_m)))\\
& = & 	g(\psi_m(0),\psi_m(-d_2(L\psi_m)),\ldots,\psi_m(-d_{{\bf k}}(L\psi_m)))
\end{eqnarray*}
is contained in the set $\{g(v)\in\R^n:v\in V\,\,\mbox{and}\,\,(v_1,\ldots,v_n)^{tr}\in B_0\}$, which is bounded due to the hypothesis. Analogously the set  
$\{f(\tilde{\psi}_m)\in\R^n:m\in\N\}$ is bounded.
Proceed as in Part 1.1.1, beginning with applications of the Bolzano-Weierstra{\ss} Theorem. 

\medskip

1.1.3. In case $U_J=U$ we get $\psi_m\in U^J\subset U=U_J$ for every $m\in\N$, and Proposition 4.2 (iv) yields $\chi_m=\psi_m-Y_J(L\psi_m)\cdot f(\psi_m)\in U_J$ and $f(\chi_m)=f(\psi_m)$ for every $m\in\N$. Proposition 4.2 (iv)
applied to $\phi\in U_J$ gives $\chi=R^J(\phi)=\phi-Y_J(L\phi)\cdot\phi'(0)\in U_J$. Now continuity yields
$f(\chi_m)\to f(\chi)$ as $m\to\infty$, and it follows that
$$
\psi_m=\chi_m+Y_J(L\chi_m)\cdot f(\chi_m)\to\chi+Y_J(L\chi)\cdot f(\chi)
$$
as $m\to\infty$.
Analogously, $\tilde{\psi}_m\to\chi+Y_J(L\chi)\cdot f(\chi)$ as $m\to\infty$.
Set $m_p=p$ for all $p\in\N$ and $x=\tilde{x}=f(\chi)\in\R^n$.

\medskip

1.2. For $\psi=\chi+Y_J(L\chi)\cdot x$ we show $\psi\in X_{fJ}$ and $R^J(\psi)=\chi$. 

\medskip

1.2.1.  Proof of $\psi\in U_J$. By Eqs. (14) and (16), $\psi(0)=\chi(0)=\phi(0)$ and $L\psi=L\chi=L\phi\in W$.
From $\phi\in X_{fJ}\subset U_J$ we obtain $L\phi\in W_J$, so for every $j\in J$, $d_j(L\psi)=d_j(L\phi)=0$, hence
$$
\psi(-d_j(L\psi))=\psi(0)=\phi(0)=\phi(-d_j(L\phi)).
$$
For $k\in K\setminus J$ we use $L\chi=L\phi\in W_J\subset W^J$ and apply Eq. (17) to $w=L\chi\in W^J$. This yields 
$$
(Y_j(L\chi)\cdot x)(-d_k(L\chi))=0,
$$
hence
$$
\psi(-d_k(L\chi))=\chi(-d_k(L\chi)).
$$
Similarly, by $R^J(\phi)=\chi$ and Eq. (17),
$$
\phi(-d_k(L\phi))=\chi(-d_k(L\phi)).
$$
It follows that 
$$
\psi(-d_k(L\psi))=\psi(-d_k(L\chi))=\chi(-d_k(L\chi))=\chi(-d_k(L\phi))=\phi(-d_k(L\phi)).
$$
Altogether, $\widehat{\psi}=\widehat{\phi}\in V$, and
thereby, $\psi\in U$. As $L\psi=L\phi\in W_J$ we also obtain $\psi\in U_J$. 

\medskip	

1.2.2. Using $X_f\ni\psi_{m_p}\to\psi\in U_J$ for $p\to\infty$ and the fact that $X_f$ is closed in $U$ we infer $\psi\in X_f\cap U_J=X_{fJ}$. Using $L\psi\in W_J\subset W^J$, the equations (16) and (15), and $\chi\in X_0$ we obtain from the definition $\psi=\chi+Y_J(L\chi)\cdot x$ the relations $ L\psi=L\chi$ and $\psi'(0)=\chi'(0)+x=x$. Hence $R^J(\psi)=\psi-Y_J(L\psi)\cdot\psi'(0)=\psi-Y_J(L\chi)\cdot x=\chi$.

\medskip

1.3. Set $\tilde{\psi}=\chi+Y_J(L\chi)\cdot\tilde{x}$. Part 1.2 with $\tilde{x}$ in place of $x$ yields $\tilde{\psi}\in X_{fJ}$ and $R^J(\tilde{\psi})=\chi$. By Proposition 4.4 the restriction of $R^J$ to $X_{fJ}$ is injective, hence
$\psi=\tilde{\psi}=\phi$. For $p\in\N$ sufficiently large we obtain $\psi_{m_p}\in N_{\phi}$ and $\tilde{\psi}_{m_p}\in N_{\phi}$, and arrive at a contradiction to Corollary 4.5 which guarantees the injectivity of $R^J$ on $N_{\phi}$. 

\medskip

2. Proof of assertion (ii). Let $\phi$ and $\psi$ in 
$$
X_f\cap(\cup_{\chi\in R^J(X_{fJ})}(R^J)^{-1}(V_{\chi}))
$$
be given with $R^J(\phi)=R^J(\psi)$. $R^J(\phi)$ is contained in $V_{\chi}$ for some $\chi\in R^J(X_{fJ})$, and $R^J(\psi)=R^J(\phi)$ is contained in the same set $V_{\chi}$. Or, both $\phi$ and  $\psi$ belong to
$X_f\cap (R^J)^{-1}(V_{\chi})$. Assertion (i)  yields $\phi=\psi$. $\Box$

\medskip

The next result establishes that $R^J$ defines a manifold chart for the solution manifold $X_f$, on a domain which contains the set $X_{fJ}$ and is an almost graph over $X_0$. This is parallel to a part of the assertions of Theorem 3.5, and the proof of Theorem 4.7 below follows its counterpart in Section 3, up to changes due to the fact that unlike $R^K$ in Section 3 the map $R^J$ for $J\neq K$ with $U_J\neq\emptyset$ can not be assumed to be the restriction of some projection onto $X_0$. The example from \cite[Section 3]{W6} 
shows that this complication is not a deficiency caused by the construction of $R^J$, and that we can not expect a neighbourhood of $X_{fJ}$ in $X_f$ to be a graph over $X_0$. In case of the example we have, in present terminology, $U_{\emptyset}=U$, and $X_{f\emptyset}=X_f\cap U_{\emptyset}=X_f$ has no graph representation with respect to any direct sum decomposition of $C^1_n$ into closed subspaces. This means in particular that a necessary condition for the injective map 
$$
X_{f\emptyset}\ni\phi\mapsto R^{\emptyset}(\phi)\in X_0
$$ 
being a 
restricted projection is violated.

\medskip

For each $\phi\in X_{fJ}$ we choose an open neighbourhood $N_{\phi}$ in $X_f$ according to Corollary 4.5, and an open neighbourhood $V_{\chi}$ of $\chi=R^J(\phi)$ in $X_0$ according to Proposition 4.6. We set
$$
N_J=\cup_{\phi\in X_{fJ}}(N_{\phi}\cap(R^J)^{-1}(V_{R^J(\phi)})).
$$

\begin{thm}
Let $J\subset K$ with $J\neq K$ and $U_J\neq\emptyset$ be given and suppose that one of the hypotheses (b), (d1b), (J) is satisfied. Then the map $R^J$ defines a diffeomorphism $R_J$ from the open neighbourhood $N_J$ of $X_{fJ}$ in $X_f$ onto the open subset $R^J(N_J)$
of the space $X_0$. We have
$$
(R_J)^{-1}(\chi)=\chi+Y_J(L\chi)\cdot f((R_J)^{-1}(\chi))\quad\mbox{for every}\quad\chi\in R^J(N_J)
$$
and the set 
$$
N_J=\{\chi+Y_J(L\chi)\cdot f((R_J)^{-1}(\chi)):\chi\in R^J(N_J)\}
$$ 
is an almost graph over $X_0$.
\end{thm}

{\bf Proof.} 1.  The set $N_J$ is open in $X_f$ and contains $X_{fJ}$. Due to Proposition 4.6 (ii) the restriction of $R^J$ to $N_J$ is injective.
Corollary 4.5 yields that locally this restriction is given by diffeomorphisms onto open subsets of the space $X_0$. It follows easily  that $R^J$ defines a diffeomorphism $R_J$ from $N_J$ onto the open subset $R^J(N_J)$ of $X_0$.

\medskip

2. Computation of $(R_J)^{-1}$. Let $\chi\in R^J(N_J)$ be given and set $\phi=(R_J)^{-1}(\chi)\in N_J\subset X_f$. Then $\chi=R^J(\phi)=\phi-Y_J(L\phi)\cdot\phi'(0)=\phi-Y_J(L\phi)\cdot f(\phi)$. Using this and Eq. (16) we get $L\phi=L\chi$, and thereby
$$ 
(R_J)^{-1}(\chi)=\phi=\chi+Y_J(L\chi)\cdot f((R_J)^{-1}(\chi)).
$$

\medskip

3. The assertion about the representation of $N_J=(R_J)^{-1}(R^J(N_J))$ is obvious from Part 2. In order to show that $N_J$ is an almost graph over $X_0$ we prove
$$
Y_J(L\chi)\cdot f((R_J)^{-1}(\chi))=0\quad\mbox{on}\quad R^J(N_J)\cap N_J
$$ 
and
$$	
Y_J(L\chi)\cdot f((R_J)^{-1}(\chi))\in C^1_n\setminus X_0	\quad\mbox{on}\quad R^J(N_J)\setminus N_J.
$$ 
	
3.1. For $\chi\in R^J(N_J)\cap N_J$, $R^J(\chi)=\chi$ (since $\chi\in X_0\cap U^J$). Using this and $\chi\in N_J$ we get $(R_J)^{-1}(\chi)=\chi$. It follows that
\begin{eqnarray*}
f((R_J)^{-1}(\chi)) & = & f(\chi)=\chi'(0)	\quad\mbox{(as}\quad\chi\in N_J\subset X_f)\\
& = & 0	\quad\mbox{(as}\quad\chi\in R^J(N_J)\subset X_0),
\end{eqnarray*}
hence $Y_J(L\chi)\cdot f((R_J)^{-1}(\chi))=0$.

\medskip

3.2. For $\chi\in R^J(N_J)\setminus N_J$ set $\phi=(R_J)^{-1}(\chi)$ and assume
$Y_J(L\chi)\cdot f(\phi)\in X_0$. By Proposition 4.2 (iii),
$0=f(\phi)$. As $\phi\in N_J\subset X_f$, $\phi'(0)=f(\phi)=0$, or $\phi\in X_0$, which yields $\phi=R^J(\phi)=\chi$. Hence $\chi\in N_J$, in contradiction to the choice of $\chi$. It follows that
$Y_J(L\chi)\cdot f((R_J)^{-1}(\chi))\in C^1_n\setminus X_0$. $\Box$

\medskip

Finally we obtain the desired almost graph diffeomorphism. This is a bit more involved than the corresponding part of the proof of Theorem 3.5, again due to the fact that $R^J$ can not be assumed to be a restriction of a projection.

\begin{thm}
Let $J\subset K$ with $J\neq K$ and $U_J\neq\emptyset$ be given and suppose that one of the hypotheses (b), (d1b), (J) is satisfied. Then then are an open subset ${\mathcal O}_J$ of the space $C^1_n$ with
$X_{fJ}\subset X_f\cap{\mathcal O}_J$ and a diffeomorphism $A_J$ from ${\mathcal O}_J$ onto an open subset of $C^1_n$ with 
\begin{eqnarray*}
A_J(\phi) & = & R^J(\phi)\quad\mbox{for every}\quad\phi\in X_f\cap{\mathcal O}_J\subset N_J,\\
A_J(\chi) & = & \chi\quad\mbox{for every}\quad\chi\in X_0\cap(X_f\cap{\mathcal O}_J),\\
& & \mbox{and}\\
A_J(X_f\cap{\mathcal O}_J) & = & X_0\cap A_J({\mathcal O}_J).
\end{eqnarray*}
\end{thm}

{\bf Proof.} 1. The set ${\mathcal O}_{(J)}=(R^J)^{-1}(R^J(N_J))$ is an open subset of the domain $U^J\subset C^1_n$. Obviously, $N_J\subset {\mathcal O}_{(J)}$ and $R^J({\mathcal O}_{(J)})\subset R^J(N_J)$. So for $\phi\in{\mathcal O}_{(J)}$ the element $R^J(\phi)$ is in the domain $R^J(N_J)$ of $(R_J)^{-1}$. We see in particular that the equation 
$$
A_{(J)}(\phi)=\phi-Y_J(L\phi)\cdot f((R_J)^{-1}(R^J(\phi)))
$$
defines a continuously differentiable map $A_{(J)}:{\mathcal O}_{(J)}\to C^1_n$.

\medskip

2. Proof of $X_f\cap{\mathcal O}_{(J)}=N_J$. The inclusion $N_J\subset X_f\cap{\mathcal O}_{(J)}$ is obvious from $N_J\subset X_f$ and 
$$
N_J\subset(R^J)^{-1}(R^J(N_J))={\mathcal O}_{(J)}.
$$ 
In order to show the reverse inclusion let $\phi\in X_f\cap{\mathcal O}_{(J)}$ be given. Then 
$$
R^J(\phi)\in R^J(N_J),
$$
or, $R^J(\phi)=R^J(\psi)$ for some
$$
\psi\in N_J\subset X_f\cap(\cup_{\chi\in R^J(X_{fJ})}(R^J)^{-1}(V_{\chi})).
$$
Hence 
$$
R^J(\phi)=R^J(\psi)\in V_{\eta}\quad\mbox{for some}\quad\eta\in R^J(X_{fJ}),
$$
and both $\phi$ and $\psi$ belong to the set $X_f\cap (R^J)^{-1}(V_{\eta})$. Proposition 4.6 (i) yields  $\phi=\psi\in N_J$.

\medskip

3. For $\phi\in N_J=X_f\cap{\mathcal O}_{(J)}$ we have
\begin{eqnarray*}
	A_{(J)}(\phi) & = & \phi-Y_J(L\phi)\cdot f((R_J)^{-1}(R^J(\phi)))=\phi-Y_J(L\phi)\cdot f(\phi)\\
	& = & \phi-Y_J(L\phi)\cdot\phi'(0)=R^J\phi.
\end{eqnarray*}

4. The set ${\mathcal O}_J=(A_{(J)})^{-1}({\mathcal O}_{(J)})$ is contained in the domain ${\mathcal O}_{(J)}$ of $A_{(J)}$. We prove $X_{fJ}\subset{\mathcal O}_J$. Let $\phi\in X_{fJ}\subset N_J$ be given. By Part 3, $A_{(J)}(\phi)=R^J(\phi)$. In order to obtain $\phi\in{\mathcal O}_J$ we must show $A_{(J)}(\phi)\in{\mathcal O}_{(J)}$, or equivalently, $A_{(J)}(\phi)\in U^J$ and $R^J(A_{(J)}(\phi))\in R^J(N_J)$.
Proof of the preceding statement: Using the definition of $R^J$ in combination with $\phi\in X_{fJ}\subset U_J$ and Proposition 4.2 (iv) we observe $R^J(\phi)\in U_J\subset U^J$. Hence $A_{(J)}(\phi)=R^J(\phi)\in U^J$, and
\begin{eqnarray*}
R^J(A_{(J)}(\phi)) & = & R^J(R^J(\phi))=R^J(\phi)-Y_J(LR^J(\phi))\cdot(R^J(\phi))'(0)\\
& = & R^J(\phi)\quad\mbox{(since}\quad (R^J(\phi))'(0)=0\quad\mbox{due to}\quad R^J(\phi)\in X_0)\\	
& \in & R^J(N_J).	
\end{eqnarray*}

5. As ${\mathcal O}_J$ is an open subset of $C^1_n$ the set $X_f\cap{\mathcal O}_J\supset X_{fJ}$ is an open neighbourhood of $X_{fJ}$ in $X_f$ which is contained in $X_f\cap {\mathcal O}_{(J)}=N_J$.

\medskip

6. For each $\rho\in{\mathcal O}_{(J)}$ we have $R^J(\rho)\in R^J(N_J)$, that is, $R^J(\rho)$ belongs to the domain of $(R_J)^{-1}$. Therefore the equation
$$
B_J(\rho)=\rho+Y_J(L\rho)\cdot f((R_J)^{-1}(R^J(\rho)))
$$
defines a continuously differentiable map $B_J:{\mathcal O}_{(J)}\to C^1_n$.

\medskip

For $\phi\in{\mathcal O}_J$ we have $A_{(J)}(\phi)\in{\mathcal O}_{(J)}$.

\medskip

Proof of $B_J(A_{(J)}(\phi))=\phi$ for every $\phi\in{\mathcal O}_J$. Let $\phi\in{\mathcal O}_J$ be given and set $\rho=A_{(J)}(\phi)\in{\mathcal O}_{(J)}$. Using the definition of $A_{(J)}$ and Eqs. (16) and (15) we get
$L\rho=L\phi$ and
$$
\rho'(0)=\phi'(0)-f((R_J)^{-1}(R^J(\phi))).
$$
Consequently,
\begin{eqnarray}
R^J(\rho) & = & \rho-Y_J(L\rho)\cdot\rho'(0)=\rho-Y_J(L\phi)\cdot\rho'(0)\nonumber\\
& = & [\phi-Y_J(L\phi)\cdot f((R_J)^{-1}(R^J(\phi)))]\nonumber\\
& & -Y_J(L\phi)\cdot[\phi'(0)-f((R_J)^{-1}(R^J(\phi)))]\nonumber\\
& = & \phi-Y_J(L\phi)\cdot\phi'(0)=R^J(\phi).\nonumber
\end{eqnarray}
It follows that
\begin{eqnarray*}
B_J(A_{(J)}(\phi)) & =	& B_J(\rho)=\rho+Y_J(L\rho)\cdot
f((R_J)^{-1}(R^J(\rho)))\\
& = & A_{(J)}(\phi)+Y_J(L\rho)\cdot
f((R_J)^{-1}(R^J(\rho)))\\
& = & [\phi-Y_J(L\phi)\cdot f((R_J)^{-1}(R^J(\phi)))]\\
& & + Y_J(L\rho)\cdot
f((R_J)^{-1}(R^J(\rho)))\\
& = & [\phi-Y_J(L\phi)\cdot f((R_J)^{-1}(R^J(\phi)))]\\
& & +Y_J(L\phi)\cdot f((R_J)^{-1}(R^J(\phi)))=\phi.
\end{eqnarray*}

\medskip

7. Part 6 yields that $A_{(J)}({\mathcal O}_J)$ is contained in the open  subset ${\mathcal O}_{JB}=(B_J)^{-1}({\mathcal O}_J)$ of the space $C^1_n$. As in Part 6 one shows
$$
A_{(J)}(B_J(\rho))=\rho\quad\mbox{for every}\quad\rho\in{\mathcal O}_{JB}.
$$
It follows that $A_{(J)}$ defines a diffeomorphism $A_J$ from ${\mathcal O}_J$ onto ${\mathcal O}_{JB}$.

\medskip

8. We have
\begin{eqnarray*}
A_J(X_f\cap{\mathcal O}_J) & \subset & A_{(J)}(X_f\cap{\mathcal O}_{(J)})\cap A_J({\mathcal O}_J)\\
& = & R^J(X_f\cap{\mathcal O}_{(J)})\cap A_J({\mathcal O}_J)\quad\mbox{(see Part 3)}\\
& \subset & X_0\cap A_J({\mathcal O}_J).
\end{eqnarray*}
In order to show the reverse inclusion let $\chi\in X_0$ with $\chi=A_J(\phi)$ for some $\phi\in{\mathcal O}_J$ be given. By the definition of ${\mathcal O}_J$, $$
\chi=A_J(\phi)=A_{(J)}(\phi)\in{\mathcal O}_{(J)}\subset U^J.
$$ 
Using $\chi\in X_0$, or $\chi'(0)=0$, we get
\begin{eqnarray*}
\chi & = & R^J(\chi)\\
& = & R^J(A_J(\phi))\in R^J({\mathcal O}_{(J)})\\
& \subset & R^J(N_J)\quad\mbox{(due to the definition of}\quad {\mathcal O}_{(J)}).
\end{eqnarray*} 
Hence $\chi=R^J(\psi)$ for some $\psi\in N_J=X_f\cap{\mathcal O}_{(J)}$. Using Part 3 we infer
$A_{(J)}(\psi)=R^J(\psi)=\chi\in{\mathcal O}_{(J)}$. Or,
$$
\psi\in (A_{(J)})^{-1}({\mathcal O}_{(J)})={\mathcal O}_J.
$$
It follows that $A_J(\psi)=A_{(J)}(\psi)=\chi$ with $\psi\in X_f\cap{\mathcal O}_J$, hence $\chi\in A_J(X_f\cap{\mathcal O}_J)$..

\medskip

9. On $X_f\cap{\mathcal O}_J\subset X_f\cap{\mathcal O}_{(J)}$
we have $A_J(\phi)=A_{(J)}(\phi)=R^J(\phi)$, see Part 3. For $\chi\in X_0\cap(X_f\cap{\mathcal O}_J)$ we infer $A_J(\chi)=R^J(\chi)=\chi$. $\Box$

\end{document}